\documentclass[pdflatex,sn-mathphys-ay]{sn-jnl}

\usepackage{graphics}
\usepackage{epsfig}
\usepackage{latexsym}
\usepackage{amsmath}
\usepackage{amssymb}
\usepackage{rotating}
\usepackage{amstext}
\usepackage{mathtools}
\usepackage{array} 
\usepackage{algorithm}
\usepackage[noend]{algpseudocode} 
\usepackage{tabularx} 
\usepackage{booktabs} 
\usepackage{dsfont}
\usepackage{longtable} 
\usepackage{enumerate} 
\usepackage{ulem}
\usepackage{tikz}
\usetikzlibrary{fit,positioning,arrows.meta}
\usepackage{pgfplots}
\pgfplotsset{compat=1.18} 
\usepackage{multicol}
\usepackage{multirow}
\usepackage{bigstrut}
\usepackage{url}
\usepackage{hyperref} 
\usepackage{placeins} 
\newcommand{\rev}[1]{#1}
\newcommand{\new}[1]{#1}

\hypersetup{citecolor=black,linkcolor=black,urlcolor=black,filecolor=black}

\makeatletter
\newcommand*{\addFileDependency}[1]{
  \typeout{(#1)}
  \@addtofilelist{#1}
  \IfFileExists{#1}{}{\typeout{No file #1.}}
}
\makeatother

\theoremstyle{thmstyleone}%
\theoremstyle{thmstyletwo}%

\theoremstyle{thmstylethree}%

\definecolor{verdeagua}{rgb}{0.0, 0.8, 0.8}

\begin{document}

\title[AI-Assisted Decision Support for Formula 1]{AI-Assisted Decision Support for Formula 1 Pit-Stop Strategy: Deep Reinforcement Learning and Mathematical Programming}

\author[1]{\fnm{Felipe} \sur{Méndez}}\email{felipe.mendez.p@ug.uchile.cl}

\author*[1]{\fnm{Charles} \sur{Thraves}}\email{cthraves@dii.uchile.cl}

\affil*[1]{\orgdiv{Department of Industrial Engineering}, \orgname{University of Chile}, \orgaddress{\street{Beauchef 851}, \city{Santiago}, \country{Chile}}}

\abstract{Pit-stop strategy in Formula~1 is a sequential decision problem under uncertainty, in which strategists decide at each lap whether to pit and which tire compound to fit. We develop a data-driven simulator, calibrated with historical race data, that captures tire degradation, drivers' recent performance, overtaking, yellow flag events, and rivals' pit-stop decisions. We train agents with five Deep Reinforcement Learning (DRL) algorithms and compare them with a Quadratic Mixed-Integer Programming (QMIP) formulation that is optimal in deterministic single-driver settings. Without uncertainty or competition, all DRL methods except TRPO achieve near-optimal performance. With a full grid of competitors and stochastic race events, the ranking of DRL methods depends on the evaluation setting: Proximal Policy Optimization (PPO) and its recurrent variant are the only methods that improve on the QMIP both in fully simulated races and when the actual 2024 races are replayed. Averaged over eight circuits, PPO improves on the QMIP by 0.33 and 0.49 positions in the two settings and Recurrent PPO by less, and in the replay PPO finishes 0.66 positions ahead of the drivers' actual strategies simulated in the same races. The QMIP performs similarly to the other methods, even though it cannot anticipate yellow flag events. Finally, we use the trained PPO policy as a look-ahead decision-support tool that shows the strategist the distribution of race outcomes under each candidate pit-stop action. Applied retrospectively to eight 2024 races, the look-ahead identifies the laps where the actual decision, the recommendation, or a third action was preferable.}

\keywords{Decision Support, Human--AI Collaboration, Deep Reinforcement Learning, Pit-Stop Strategy, Mixed-Integer Programming, OR in Sports}

\maketitle

\section{Introduction}\label{sec:intro}

Many operational systems require high-stakes sequential decisions to be made under uncertainty, strategic interaction, and severe time pressure. In such settings, human decision-makers increasingly rely on AI-based decision-support systems, but fully autonomous policies are often difficult to trust, interpret, and validate. This paper studies this challenge in the context of Formula~1 pit-stop strategy, where teams must decide under tire degradation, traffic, competitors' actions, and stochastic race interruptions. 

Although Formula~1 is a motorsport application, the underlying decision structure is shared by many operational settings in which human decision-makers must act under uncertainty, strategic interaction, and limited time, such as dynamic vehicle dispatching \citep{pillac2013review}, disruption management in transportation networks \citep{clausen2010disruption}, and real-time resource allocation \citep{powell2016perspectives}.

Pit-stop strategy is a key determinant of performance in Formula~1 racing. At each lap, a driver must decide whether to remain on track or perform a pit stop, trading off the immediate time loss associated with pitting against the potential performance gains from fresher tires. This decision is further complicated by the availability of different tire compounds, namely soft, medium, and hard, which exhibit distinct trade-offs between speed and degradation. Softer compounds provide faster lap times but deteriorate more rapidly, while harder compounds are more durable but slower. In addition, race regulations require drivers to use at least two different tire compounds under dry conditions, further constraining feasible strategies.

Stochastic race events also affect these decisions. In particular, yellow-flag periods, including Safety Car (SC) and Virtual Safety Car (VSC) events, can substantially alter the relative cost of pit stops. During these periods, cars are forced to reduce speed and overtaking on track is restricted, which reduces the relative time lost in the pit lane and can create strategic opportunities. Effective decisions therefore depend on the current race state and on the possibility of future disruptions.

Incorporating competition against rival cars further complicates the problem through both physical interactions and strategic behavior. From a physical standpoint, traffic conditions affect lap times, as overtaking is often difficult and inherently stochastic, even when relative pace differences are known. From a strategic perspective, each driver’s decision influences and is influenced by the actions of others, leading to a high-dimensional dynamic game. With approximately 20 drivers making pit-stop and tire-compound decisions over the course of a race, the joint decision space grows rapidly with the number of laps and drivers, making exact modeling approaches computationally challenging.

The combination of stochastic race dynamics, operational constraints, and strategic interactions makes pit-stop optimization difficult. This difficulty is particularly relevant for AI-assisted decision support: exact optimization can provide high-quality solutions and structural insight, but becomes computationally impractical in realistic multi-driver settings, while learning-based policies can generate fast recommendations once trained but require careful validation and interpretation. This paper addresses this trade-off by studying deep reinforcement learning (DRL) methods for Formula~1 pit-stop strategy and evaluating them against an exact mathematical programming benchmark. From this perspective, learning-based policies are not intended to replace human strategists, but to support them with fast recommendations and interpretable summaries of alternative actions.

Existing approaches to pit-stop strategy include simulation-based methods, supervised learning models trained on historical race decisions, exact optimization approaches based on dynamic programming or mathematical programming, and, more recently, reinforcement learning methods. Simulation and supervised learning can capture realistic race features and generate fast recommendations, but they do not necessarily identify optimal strategies. Exact optimization methods provide benchmarks and structural insight, but their scalability is limited in realistic multi-driver settings. Learning-based policies scale to high-dimensional sequential decision problems, but their use as decision-support tools raises questions of reliability and interpretability that motivate our study. Section~\ref{sec:related_literature} reviews these streams in detail.

Three challenges remain open for AI-assisted decision support in pit-stop strategy. First, learning-based policies should be validated against exact optimization benchmarks whenever such benchmarks are available. Second, their performance under competition and stochastic race events should be checked against the realizations of actual races, which, to our knowledge, has not been done for DRL. Third, recommendations should be translated into interpretable summaries that allow human strategists to compare alternative actions and assess their associated risks.

To address these challenges, we investigate the following three research questions: 
(i) To what extent can DRL-based pit-stop policies reproduce optimal decisions in settings where exact solutions are available? 
(ii) How reliable are learning-based recommendations under competition and multiple sources of uncertainty, in simulated races and in replays of the actual 2024 races? 
(iii) How can the outcomes of alternative pit-stop actions be evaluated and visualized to help human strategists understand the risks and expected outcomes of competing race decisions?

The main contributions of this paper are as follows:
\begin{itemize}
\item We develop a data-driven race simulation environment calibrated with public historical race data from 2021 to 2024. The environment captures key sources of operational uncertainty, including lap-time variability, driver interactions, overtaking dynamics, and yellow-flag events. It also incorporates predictive models of rivals' pit-stop behavior and updates driver-specific performance estimates using only information available prior to each race, so that the whole procedure can be run before a Grand Prix.

\item We evaluate five DRL methods against an exact QMIP benchmark. In simplified deterministic settings, where optimal solutions can be computed, this comparison allows us to quantify the optimality gap of learning-based policies. In realistic full-grid scenarios with competition and multiple sources of uncertainty, we compare learning-based recommendations with an exact optimization benchmark that uses a simplified race representation. The evaluation covers eight circuits of the 2024 season, in both fully simulated races and replays of the actual races. Performance relative to the QMIP varies across DRL methods, so the choice of algorithm matters.

\item We show how the resulting policies can support human race strategists through look-ahead simulations of alternative pit-stop decisions. For selected race states, the framework evaluates and visualizes the distribution of final outcomes under each feasible action, allowing decision-makers to compare expected performance, downside risk, and outcome variability. We also apply it retrospectively to review the pit-stop decisions taken in the 2024 races.
\end{itemize}

\rev{The remainder of the paper is organized as follows. Section~\ref{sec:related_literature} reviews the relevant literature. Section~\ref{sec:environment} describes the race simulation environment and its calibration based on historical race data. Section~\ref{sec:agent} formulates the pit-stop decision problem, including the state and action spaces, rivals' behavior, and the reward function. Section~\ref{sec:QMIP} introduces an exact quadratic mixed-integer programming (QMIP) formulation for simplified instances. Section~\ref{sec:results} validates the learned policies against the QMIP and in full-grid race scenarios. Section~\ref{sec:decision_support} presents the look-ahead decision-support analysis, and Section~\ref{sec:conclusions} concludes. The appendices follow the references.}

\section{Related Literature}\label{sec:related_literature}

\subsection{Race Strategy and Pit-Stop Optimization}
Different approaches have been proposed to address the pit-stop optimization problem, ranging from exact methods to simulation-based and learning-based techniques. Exact approaches, such as mathematical programming and dynamic programming (DP), provide optimal solutions but are typically limited to simplified settings due to the exponential growth of the state space. Prior to the development of such methods, simulation was widely used to evaluate candidate strategies by modeling race dynamics that involve uncertainty, including overtaking, tire degradation, lap time variability, safety car deployments, Drag Reduction System (DRS),\footnote{The \textit{Drag Reduction System}, introduced in 2011, is a car feature that allows a driver to open the rear wing to reduce air resistance on certain straights when trailing another car within one second.} and driver interactions \citep{Bekker_Lotz, phillips2014building, heilmeier2018race, sulsters2018simulating, heilmeier2020application}. While simulation enables the incorporation of realistic and complex racing events, it is computationally intensive and requires evaluating a large number of scenarios to identify high-performing strategies.

Machine learning approaches have also been explored as decision-support tools. Supervised methods such as Random Forests, Support Vector Machines, and LASSO have been trained on historical race data to emulate pit-stop decisions \citep{tulabandhula2014tire, choo2015real, fatima2023deep}. More recently, \cite{Sasikumar2025} propose a supervised deep learning framework using Bi-LSTM and related architectures to predict pit stop timing from telemetry data, achieving strong classification performance. These approaches can provide real-time recommendations but are inherently limited by the quality of historical decisions, as they tend to replicate past strategies rather than identify optimal ones. Hybrid approaches combining optimization and learning have also been proposed; for instance, \cite{heilmeier2020virtual} integrate neural networks with optimization models to guide decision-making.

Exact solution methods based on dynamic programming (DP) have recently been proposed for pit-stop optimization problems. In a single-driver setting, \cite{heine2023optimization} formulate the problem as a stochastic DP in which, at each lap, the driver decides whether to change tires based on the current race state. Despite incorporating stochastic elements such as yellow flag events, the absence of competition allows for the practical enumeration of all relevant states. Extensions to multi-driver settings have been explored using game-theoretic formulations. In particular, \cite{Aguad_Thraves} model the interaction between two drivers as a feedback Stackelberg game, showing that the problem remains computationally tractable despite a state space on the order of billions. However, these approaches do not scale to realistic race settings involving more than two drivers, as the state space grows exponentially with the number of agents. This limitation motivates the use of approximate methods, such as DRL, which can handle large state spaces while maintaining a relatively small action space.

\subsection{Learning-Based Methods for Race Strategy}
The comparison between DRL and classical optimization approaches has attracted growing interest across multiple operations research domains. In inventory management, \cite{gijsbrechts2022deep} provide a rigorous evaluation of DRL against optimal and near-optimal analytical policies in classic intractable inventory problems, while \cite{dehaybe2024deep} evaluate DRL under non-stationary uncertain demand. In power systems scheduling, \cite{omalley2023reinforcement} conduct a head-to-head comparison of MIP and reinforcement learning for unit commitment, finding that MIP consistently produces lower-cost schedules, though RL's key advantage lies in generating solutions practically instantly regardless of problem size. These results suggest that the relative performance of DRL and classical optimization methods is highly problem-dependent. However, comparable evaluations against exact optimization benchmarks remain limited in the context of pit-stop strategy optimization.

DRL has recently been applied to high-dimensional sequential decision problems. In motorsports, DRL has been applied to trajectory optimization, energy management, and refueling strategies \citep{remonda2021formula, hoeppke2022optimal, liu2023learn, evans2023safe, liu2021formula, boettinger2023mastering}. 

In the context of pit-stop strategies, \cite{piccinotti2021online} simulate lap times calibrated with XGBoost and evaluate different implementations of Monte Carlo Tree Search (MCTS) to identify optimal sequences of actions. Among the methods considered, an open-loop Q-learning approach achieves the best performance. More recent work has focused on DRL methods. Using Deep Q-Networks (DQN), \cite{marinaro2023smarl} compares single-agent, multi-agent, and hybrid learning approaches, showing that the hybrid setting achieves strong performance with reduced data requirements. \cite{rebe2024toma} evaluates DQN, QR-DQN, and A2C under different reward specifications based on positions and final outcomes, finding that A2C consistently outperforms the alternatives. Similarly, \cite{thomas2025explainable} propose a recurrent DRL architecture that integrates simulation and live-race data, using reward functions based primarily on championship points. In the first author's Master's thesis \citep{mendez2025thesis}, which this paper extends with a revised simulator, five DRL algorithms were compared with the QMIP in the single-driver deterministic case. \cite{fieni2026benchmarking}, in work independent of ours, formulate a single-driver problem at the Bahrain circuit that also allocates energy, solve it as a mixed-integer nonlinear program, and compare the solution with a Soft Actor-Critic agent. Their setting has no competitors and no stochastic race events. More recently, \cite{fieni2026learning} propose a multi-agent reinforcement learning framework that jointly optimizes energy management, pit stop timing, and tire selection using a self-play training scheme, demonstrating that agents adapt their strategies in response to competitors' actions. However, the framework is evaluated in pairwise duels between its own agents at a single circuit, without stochastic race events. 

\subsection{AI-Assisted Decision Support and Learning-Based Optimization}
The use of learning-based policies as decision-support tools connects this work to the broader literature on human--AI decision-making. In AI-assisted decisions, a recommendation should also help the decision-maker judge when to follow, question, or override it. A recurring theme in this literature is appropriate reliance: users should defer to correct recommendations and rely on their own judgment when the model is likely to err \citep{schemmer2022should}, and reliance and decision quality are distinct outcomes that explanations and additional decision-support information can help align \citep{schoeffer2025ai, leichtmann2023effects}. A complementary line of work frames decision support as the visualization and exploration of trade-offs and outcome uncertainty across a set of candidate solutions, rather than the delivery of a single recommendation; in the context of multi-objective optimization, \cite{osika2023lies} survey such methods and identify visualization, uncertainty exploration, and interactivity as core mechanisms for helping a decision-maker choose among alternatives. These ideas have also been instantiated in real-time, safety-critical operations: \cite{abbas2025analyzing} integrate deep reinforcement learning into a human-in-the-loop decision-support system for industrial control rooms, providing operators with intervention strategies adapted to the state of both the process and the operator under information overload. Pit-stop strategy shares the defining features of these settings: recommendations must be interpreted under severe time pressure, alongside information not fully encoded in the simulator, and with limited opportunity to revise a decision once made. Our look-ahead visualization builds directly on this perspective, presenting the strategist not with a single prescribed action but with the full distribution of outcomes under each feasible alternative, so that expected performance and downside risk can be weighed jointly.

Taken together, this body of work shows that learning-based methods are promising for race strategy, yet no existing study brings together the three elements needed for reliable decision support: benchmarking against an exact optimum, evaluation against actual race realizations, and interpretable comparison of the actions available at each lap.
\section{Race Simulation Environment}\label{sec:environment}

We begin by describing the problem setting. Consider a given driver racing in a Formula~1 Grand Prix at a specific circuit. Let $N$ denote the number of laps, $\mathcal{D}$ the set of drivers, $\mathcal{C}=\{\text{C0},\text{C1},\text{C2},\text{C3},\text{C4},\text{C5}\}$ the set of tire compounds\footnote{For each GP, a set of three compounds is selected from these.}, and $\mathcal{Y}$ the set of past seasons for which lap-time data are available for each driver, Grand Prix, and lap, as well as from previous races in the current season. The objective is to construct a simulation environment that generates realistic lap times and captures race dynamics, enabling the optimization of pit-stop strategies via DRL.

To train the reinforcement learning agent, it is essential to accurately model lap times, as they determine both the state transitions and the reward signal. To this end, we calibrate a regression model for each circuit using historical lap-time data. We consider two scenarios: one under clear track conditions and another under VSC or SC periods, to account for the differences in lap times across these settings. In addition, we incorporate data from previous races in the same season to update driver-specific effects and capture recent performance. Finally, we use a logit model to simulate overtaking events.

We use race data from the Python library \textit{FastF1} \citep{fastf1}, which provides detailed information for the 2021--2024 seasons. We use \cite{openf1} for tire data that are missing from or inconsistent with FastF1. The data and code supporting this study are publicly available at \cite{mendez_thraves_2025_data}. The dataset includes lap times, telemetry, car positions, tire usage, weather conditions, event schedules, and session results. 
\new{In the lap-time models, we exclude race editions in which more than 10\% of the laps were run on intermediate or wet tires, as well as laps on these tires in the remaining editions. The overtaking and rivals' models exclude only the laps on these tires. As a result, the Emilia Romagna GP, whose 2021 and 2022 editions were wet and whose 2023 edition was cancelled, is not considered.}
\new{Tire construction changed with the 2022 technical regulations, so a given tire age in 2021 does not correspond to the same wear as in later seasons. In all models that use tire age, we multiply the 2021 tire age by a factor $f=0.78$, common to all compounds and estimated from lap times and pit-stop decisions (see Appendix~\ref{app:tirelife2021}).} Tire age does not count laps under SC, VSC or red flag, and the age of a set used earlier in the weekend counts only its push laps, within 10\% of the driver's best lap of the session.
In addition, we construct variables capturing the time gaps between consecutive cars, which are critical for modeling race dynamics, as both following closely behind another car and having a car directly ahead can lead to time losses.

\subsection{Lap-Time Model under Clear-Track Conditions (RegLTC)}\label{sec:cleartrack}

\subsubsection{Base model:}\label{sec:cleartrack_base}

We use a linear regression model to predict lap times using historical data from previous races at the same circuit. A separate model is estimated for each circuit, as circuits exhibit substantial heterogeneity, particularly in how tire degradation affects lap-time performance. The regression model is specified as follows:
\begin{equation}\label{eq:reg}
\begin{aligned}
y = & \ \beta_0 + \beta_1 X_{\text{Lap}} + \beta_2 X^2_{\text{Lap}} + \beta_3 X_{\text{PitIn}} + \beta_4 X_{\text{PitOut}}
+ \sum_{d \in \mathcal{D}} \beta_d X_d + \sum_{y \in \mathcal{Y}} \beta_y X_y +\beta_{DRS}X_{DRS} \\
& + \sum_{c \in \mathcal{C}} \beta_{c,1} X_c + \sum_{c \in \mathcal{C}} \beta_{c,2} X_c X_{\text{TireLife}} + \beta_{\Delta_1} e^{-X_{\Delta_1}} + \beta_{\Delta_2} e^{-X_{\Delta_2}} + \beta_{\text{1pos}} X_{\text{1pos}},
\end{aligned}
\end{equation}
where each observation corresponds to a driver–year–lap triplet for the circuit under analysis. The intercept represents the lap time of a reference driver and tire compound on fresh tires with initial fuel load. The terms $\beta_1 X_{\text{Lap}}$ and $\beta_2 X_{\text{Lap}}^2$ capture the effect of fuel consumption over the race, allowing for non-linear improvements in lap times as the car becomes lighter. Pit-related time losses are modeled through $\beta_3 X_{\text{PitIn}}$ and $\beta_4 X_{\text{PitOut}}$. Driver- and year-specific heterogeneity are captured via fixed effects $\{\beta_d\}_{d \in \mathcal{D}}$ and $\{\beta_y\}_{y \in \mathcal{Y}}$, respectively. The effect of DRS availability is represented by $\beta_{DRS} X_{DRS}$. Tire effects are modeled through a compound-specific intercept $\beta_{c,1}$ and an interaction term $\beta_{c,2} X_c X_{\text{TireLife}}$, capturing both baseline performance and degradation. As tire usage increases, lap times are expected to increase, implying positive slope coefficients, with softer compounds exhibiting faster degradation than harder ones. Conversely, softer compounds are expected to have lower intercepts due to their superior initial grip. The terms $\beta_{\Delta_1} e^{-X_{\Delta_1}}$ and $\beta_{\Delta_2} e^{-X_{\Delta_2}}$ capture the impact of nearby cars, where the exponential specification ensures a decreasing effect as the time gap increases. Finally, $\beta_{\text{1pos}} X_{\text{1pos}}$ accounts for additional time associated with the starting position in the first lap. Table~\ref{tab:variables_regression} provides a detailed description of the independent variables. \new{At some circuits the compound coefficients $\beta_{c,1}$ and $\beta_{c,2}$ are poorly identified, either because a compound was rarely used or because tire age is highly correlated with the lap number. We therefore replace them with the estimates of a model fitted jointly across circuits, in which softer compounds are faster when new and degrade faster, and each circuit's degradation level is shrunk toward the common mean (see Appendix~\ref{app:regTC_tirecoef_matrix_completion}).}
\begin{table}[htbp]
\centering
\begin{tabular}{>{\centering\arraybackslash}m{1.5cm} >{\centering\arraybackslash}m{2.2cm} >{\raggedright\arraybackslash}m{8cm}}
\toprule
Variable & Values & Description \\
\midrule
$X_{\text{Lap}}$ & $\{1, \dots, N\}$ & Lap number \\
$X_{\text{PitIn}}$ & $\{0, 1\}$ & 1 if the car is entering the pit lane in the lap, 0 otherwise \\
$X_{\text{PitOut}}$ & $\{0, 1\}$ & 1 if the car is exiting the pit lane in the lap, 0 otherwise \\
$\{X_d\}_{d \in \mathcal{D}}$ & $\{0, 1\}$ & 1 if the driver is $d$, 0 otherwise \\
$\{X_y\}_{y \in \mathcal{Y}}$ & $\{0, 1\}$ & 1 if the year is $y$, 0 otherwise \\
$X_{DRS}$ & $\{0, 1\}$ & 1 if the car ahead is within one second, 0 otherwise \\
$\{X_c\}_{c \in \mathcal{C}}$ & $\{0, 1\}$ & 1 if tire compound $c$ is being used, 0 otherwise \\
$X_{\text{TireLife}}$ & $\{1, \dots, N\}$ & Tire age (in laps) \\
$X_{\Delta_1}$ & $\mathbb{R}_{+}$ & Time gap (seconds) to the car behind \\
$X_{\Delta_2}$ & $\mathbb{R}_{+}$ & Time gap (seconds) to the car ahead \\
$X_{\text{1pos}}$ & $\{1, \dots, 20\}$ & Starting position in the first lap (0 otherwise) \\
$X_{\text{VSC}}$ & $\{0, 1\}$ & 1 if there is a VSC during the lap, 0 otherwise \\
$X_{\text{SC}}$ & $\{0, 1\}$ & 1 if there is a SC during the lap, 0 otherwise \\
\bottomrule
\end{tabular}
\caption{Description of the independent variables in the regression model}
\label{tab:variables_regression}
\end{table}

\subsubsection{Removing outliers}\label{sec:outliers}

In some cases, observed lap times may differ considerably from the values predicted by the linear regression model. This can be explained by factors not captured in the model, many of which are unobserved. For instance, a driver may make a mistake in a corner, resulting in a loss of more than one second and a large residual. To mitigate the impact of such events, we remove these observations from the dataset to prevent them from biasing the estimation of the model coefficients. This step is particularly important for the tire compound coefficients, as inaccurate estimates would distort the agent’s understanding of tire performance and may lead to unrealistic strategies.

We first discard the laps that are at least 15 seconds slower than the median lap time of their race, except for in-laps, out-laps and the first lap. We then adopt an iterative procedure to remove outliers. First, we estimate the linear regression model and identify observations whose residuals lie more than three interquartile ranges (IQR) below the first quartile or above the third quartile of the residual distribution. If any such observations are found, they are removed from the dataset and the regression is re-estimated. This process is repeated until no additional outliers are detected (at most ten iterations). The procedure is summarized in Algorithm~\ref{alg:IQR}.
\begin{algorithm}[htbp]
\caption{Outlier elimination with IQR}\label{alg:IQR}
\begin{algorithmic}[1]
\Require{$\mathcal{I} \gets \mathcal{D}\times \mathcal{Y}\times \{1, \ldots, N\}$ and covariates-laptimes $\{x^{(i)},y^{(i)}\}_{i\in\mathcal{I}}$}
\While{true}
    \State{$\hat{y}_{i} \gets \text{LinReg} \left( \left\{ \left( x^{(i)}, y^{(i)} \right) \right\}_{i \in \mathcal{I}} \right)$ using Equation~\eqref{eq:reg}}
    \State{$\hat{e}_i \gets y_i - \hat{y}_i \quad \forall i \in \mathcal{I}$}
    \State{$Q_3 \gets 75\% \text{ percentile of } \left\{ \hat{e}_i \right\}_{i \in \mathcal{I}}$}
    \State{$Q_1 \gets 25\% \text{ percentile of } \left\{ \hat{e}_i \right\}_{i \in \mathcal{I}}$}
    \State{$\text{IQR} \gets Q_3 - Q_1$}
    \State{$\mathcal{J} \gets \left\{ i \in \mathcal{I} \mid \hat{e}_i \notin \left[ Q_1 - 3 \cdot \text{IQR}, Q_3 + 3 \cdot \text{IQR} \right] \right\}$}
    \State{$\mathcal{I} \gets \mathcal{I} \setminus \mathcal{J}$}
    \If{$|\mathcal{J}| = 0$}
        \State \textbf{break}
    \EndIf
\EndWhile
\end{algorithmic}
\end{algorithm}
It is important to note that this procedure does not directly remove outliers in lap times, but instead excludes observations with large residuals. The iterative nature of the algorithm ensures that such observations are progressively identified, including those that may initially be masked by other extreme values. The primary objective of this procedure is to prevent parameter estimates from being biased, with any improvement in model fit arising as a secondary effect.

\subsubsection{Updating drivers' coefficients based on recent performance}\label{sec:update}

Up to this point, the performance of each driver at the target circuit, captured by $\beta_{d}$ in Equation~\eqref{eq:reg}, is estimated using lap-time data from previous years of the same circuit. However, driver performance varies over time, both across seasons and across races within the same season. This variability can be explained by factors such as changes in regulations, updates in car components, and drivers adjusting to their cars, among others. Figure~\ref{fig:norris_hamilton_mexico} illustrates how a driver's performance can vary across different years. Therefore, recent races in the current season provide valuable information to better capture these dynamics in the model.

\usepgfplotslibrary{groupplots}
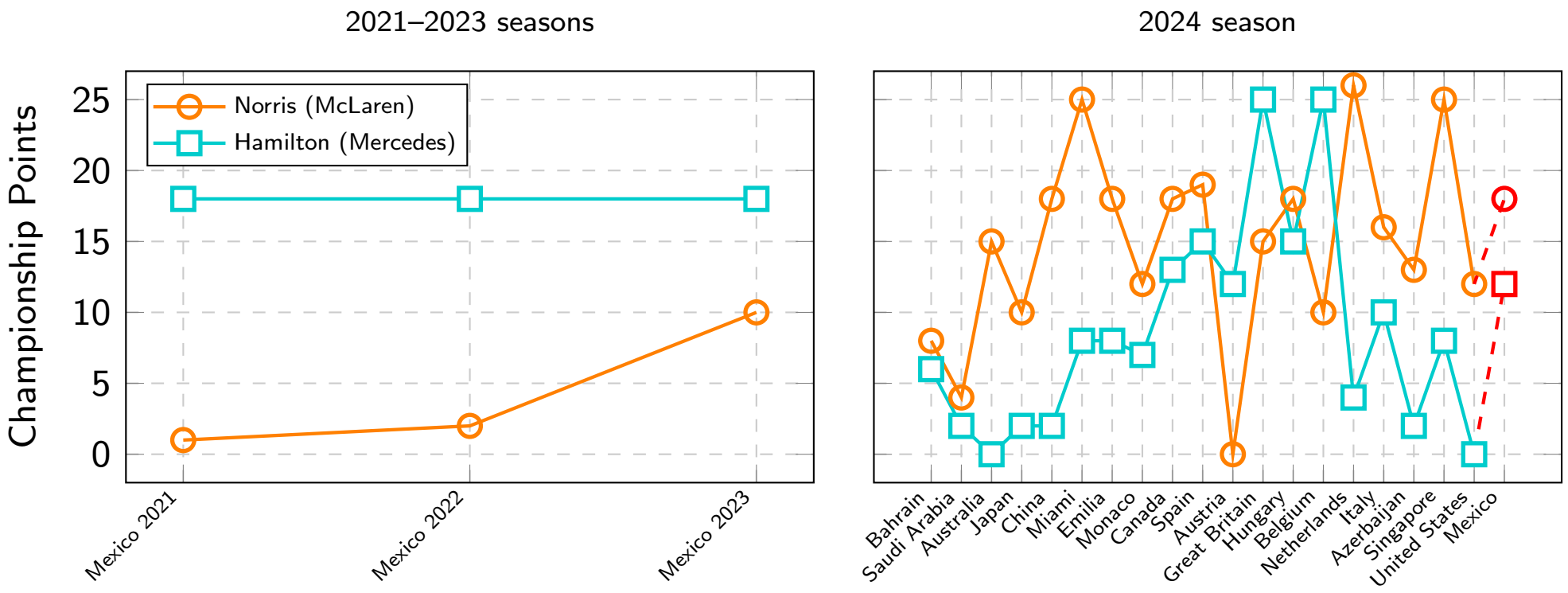
\begin{figure}[H]
\centering
\begin{tikzpicture}
\begin{groupplot}[
    group style={
        group size=2 by 1,
        horizontal sep=0.5cm,
        ylabels at=edge left,
        yticklabels at=edge left
    },
    width=7.3cm,
    height=5cm,
    ylabel={Championship Points},
    ymin=-2, ymax=27,
    ytick={0,5,10,15,20,25},
    grid=both,
    major grid style={dashed, gray!40},
    legend style={
        font=\footnotesize\sffamily,
        fill=white,
        draw=black,
        inner sep=1pt,
        nodes={scale=0.8, transform shape},
        cells={anchor=west}
    },
    title style={font=\footnotesize\sffamily},
    label style={font=\scriptsize\sffamily},
    tick label style={font=\scriptsize\sffamily},
    yticklabel={\pgfmathprintnumber[assume math mode=true]{\tick}}, 
    xticklabel style={
        rotate=45,
        anchor=east,
        font=\tiny\sffamily
    }
]

\nextgroupplot[
    symbolic x coords={
        Mexico 2021,
        Mexico 2022,
        Mexico 2023
    },
    xtick=data,
    title={2021--2023 seasons},
    legend pos=north west
]

\addplot[
    color=orange,
    mark=o,
    thick,
    mark options={
        scale=1.3,
        fill=white,
        draw=orange,
        line width=1
    }
] coordinates {
    (Mexico 2021,1)
    (Mexico 2022,2)
    (Mexico 2023,10)
};
\addlegendentry{Norris (McLaren)}

\addplot[
    color=verdeagua,
    mark=square*,
    thick,
    mark options={
        scale=1.3,
        fill=white,
        draw=verdeagua,
        line width=1
    }
] coordinates {
    (Mexico 2021,18)
    (Mexico 2022,18)
    (Mexico 2023,18)
};
\addlegendentry{Hamilton (Mercedes)}

\nextgroupplot[
    symbolic x coords={
        Bahrain, Saudi Arabia, Australia, Japan, China,
        Miami, Emilia, Monaco, Canada, Spain, Austria,
        Great Britain, Hungary, Belgium, Netherlands,
        Italy, Azerbaijan, Singapore, United States,
        Mexico
    },
    xtick={
        Bahrain, Saudi Arabia, Australia, Japan, China,
        Miami, Emilia, Monaco, Canada, Spain, Austria,
        Great Britain, Hungary, Belgium, Netherlands,
        Italy, Azerbaijan, Singapore, United States,
        Mexico
    },
    title={2024 season},
]

\addplot[
    color=orange,
    mark=o,
    thick,
    mark options={
        scale=1.3,
        fill=white,
        draw=orange,
        line width=1
    }
] coordinates {
    (Bahrain,8)
    (Saudi Arabia,4)
    (Australia,15)
    (Japan,10)
    (China,18)
    (Miami,25)
    (Emilia,18)
    (Monaco,12)
    (Canada,18)
    (Spain,19)
    (Austria,0)
    (Great Britain,15)
    (Hungary,18)
    (Belgium,10)
    (Netherlands,26)
    (Italy,16)
    (Azerbaijan,13)
    (Singapore,25)
    (United States,12)
};

\addplot[
    red,
    thick,
    dashed,
    forget plot
] coordinates {
    (United States,12)
    (Mexico,18)
};

\addplot[
    only marks,
    mark=o,
    mark options={
        scale=1.3,
        fill=white,
        draw=red,
        line width=1
    },
    forget plot
] coordinates {
    (Mexico,18)
};

\addplot[
    color=verdeagua,
    mark=square*,
    thick,
    mark options={
        scale=1.3,
        fill=white,
        draw=verdeagua,
        line width=1
    }
] coordinates {
    (Bahrain,6)
    (Saudi Arabia,2)
    (Australia,0)
    (Japan,2)
    (China,2)
    (Miami,8)
    (Emilia,8)
    (Monaco,7)
    (Canada,13)
    (Spain,15)
    (Austria,12)
    (Great Britain,25)
    (Hungary,15)
    (Belgium,25)
    (Netherlands,4)
    (Italy,10)
    (Azerbaijan,2)
    (Singapore,8)
    (United States,0)
};

\addplot[
    red,
    thick,
    dashed,
    forget plot
] coordinates {
    (United States,0)
    (Mexico,12)
};

\addplot[
    only marks,
    mark=square*,
    mark options={
        scale=1.3,
        fill=white,
        draw=red,
        line width=1
    },
    forget plot
] coordinates {
    (Mexico,12)
};

\end{groupplot}
\end{tikzpicture}
\caption{Race points of Hamilton and Norris at the Mexican Grand Prix in 2021--2023 and in the 2024 races up to that Grand Prix. Red markers show the 2024 Mexican Grand Prix.}
\label{fig:norris_hamilton_mexico}
\end{figure}
\new{To account for these dynamics, we predict each driver's coefficient for the target race from information available before it, and use these predictions in place of the coefficients $\{\beta_d\}_{d\in\mathcal{D}}$ of Equation~\eqref{eq:reg}. The prediction model is estimated with all drivers and races of 2021--2024 that precede the target race:}
\begin{equation}\label{eq:regcoef}
\new{\frac{z}{T^{\text{pole}}} = \alpha_0 + \sum_{i=1}^{3}\alpha_{\text{Prev}_i} X_{\text{Prev}_i} + \alpha_{\text{FP}} X_{\text{FP}} + \alpha_{\text{Q}} X_{\text{Q}},}
\end{equation}
\new{where each observation is a driver--race pair, $z$ is the driver's coefficient estimated with Equation~\eqref{eq:reg} using only that race, and $T^{\text{pole}}$ is the pole lap time of that race weekend. The variables $X_{\text{Prev}_i}$ are the driver's coefficients in their three most recent races, each divided by the pole lap time of its race, and $X_{\text{FP}}$ and $X_{\text{Q}}$ are the relative differences between the driver's best lap time in the free practice and qualifying sessions of the same weekend and that of the reference driver. The prediction is multiplied by the pole lap time of the target race. The treatment of missing sessions, drivers with few previous races and races with few clean laps, as well as alternative specifications, are in Appendix~\ref{app:adjust_driver}.}

The resulting linear regression used to generate lap times under clear-track conditions, denoted by \textbf{RegLTC}, is obtained by iteratively applying Equation~\eqref{eq:reg} as described in Algorithm~\ref{alg:IQR}, followed by updating the driver-specific coefficients using Equation~\eqref{eq:regcoef}.

To incorporate stochasticity, lap times are sampled from a Normal distribution with mean equal to the predicted value and variance $\hat{\sigma}^2 (1 + X_0^{T}(X^{T}X)^{-1}X_0)$, the predictive variance of the regression model, where $\hat{\sigma}^2$ is the estimated residual variance, $X_0$ is the feature vector of the new observation, and $X$ is the design matrix used for estimation. Since the driver coefficient predicted by Equation~\eqref{eq:regcoef} is also imperfect, in each simulated race it is perturbed by a draw from a zero-mean normal distribution, kept fixed during the race. Its standard deviation is that of the prediction error in the eight 2024 races of Section~\ref{sec:agent_competition}, net of the sampling error of the coefficient estimated with each race.

\subsection{Lap-Time Model under Safety Car and Virtual Safety Car (RegLTSC)}

\new{Under VSC or SC conditions, lap times are mainly determined by the imposed speed limit, and tire age and compound have little effect. We therefore model them as a multiple of the reference clear-track lap time of circuit $g$, $L^{(g)}$, and the pit-lane time loss as a fraction of the corresponding clear-track loss:}
\begin{equation}\label{eq:reg_yf}
\new{y = L^{(g)}\left(\rho_{\text{VSC}} X_{\text{VSC}} + (\rho_{\text{SC}}+u_g) X_{\text{SC}}\right) + \beta_{\text{PitIn}}' X_{\text{PitIn}} + \beta_{\text{PitOut}}' X_{\text{PitOut}},}
\end{equation}
\new{with $\beta_{\text{PitIn}}'=\kappa_{\text{in}}G^{(g)}_{\text{in}}$ and $\beta_{\text{PitOut}}'=\kappa_{\text{out}}G^{(g)}_{\text{out}}$, where $G^{(g)}_{\text{in}}$ and $G^{(g)}_{\text{out}}$ are the clear-track time losses of the in-lap and out-lap at circuit $g$, and $u_g$ is a circuit effect under SC. The model is estimated jointly over all circuits using the 2021--2023 laps run entirely under VSC or SC, so it does not require yellow flag events to be observed at every circuit (see Appendix~\ref{app:regLTYF}). Under SC, the field regroups progressively: in the second lap the gaps between consecutive cars are halved, and from the third lap the cars run in line, $1.5$ seconds apart. Pit stops under SC are only allowed in its first lap.}

The resulting regression model used to generate lap times under yellow flag conditions is denoted as \textbf{RegLTYF}.

\subsection{Overtaking Model (LogitOvertake)}\label{sec:overtake}

The components introduced in the previous sections allow us to simulate lap times. However, in practice, executing an overtaking maneuver requires more than a trailing car being sufficiently fast to close the gap with the car ahead; it must also be able to complete the pass. In other words, overtakes may not occur even when lap-time predictions suggest they should. To model overtaking behavior more realistically and consistently with observed data, we estimate the following logistic regression model:
\begin{equation}\label{eq:logitovertake}
\ln\left(\frac{p}{1-p}\right) = \gamma_{0} + \gamma_{\text{PredGap}} X_{\text{PredGap}},
\end{equation}
where $p$ denotes the probability that an overtake occurs for a given car and lap. The model is estimated per Grand Prix on pairs of consecutive cars, on laps without SC, VSC or red flag. The dependent variable is binary, taking value 1 if the trailing car finishes the lap ahead of the car it was following and 0 otherwise. The explanatory variable $X_{\text{PredGap}}$ represents the predicted time gap between two consecutive cars at the end of the lap. Specifically, it is computed as the current time gap at the beginning of the lap, plus the predicted lap time of the trailing car, minus the predicted lap time of the leading car. Predicted lap times are obtained from \textbf{RegLTC}.

It is worth noting that Equation~\eqref{eq:logitovertake} does not guarantee an overtake even when the predicted gap is negative (i.e., when the trailing car would theoretically pass the leading car), since $p<1$. To ensure physically consistent simulations, we impose that, following an overtake, the resulting time gap between consecutive cars is at least $t_{\text{min\_gap}}$.

\rev{Finally, overtaking on track is not permitted under SC or VSC conditions. Accordingly, the LogitOvertake model is applied only under clear-track conditions; during yellow flag periods, position changes occur solely as a consequence of pit stops and are resolved deterministically from the resulting time gaps.}

This model is referred to as \textbf{LogitOvertake}.

\section{Pit-Stop Decision Model}\label{sec:agent}

Near the end of each lap, the agent decides, based on the current race state, whether to stay on track or perform a pit stop, selecting one of the available tire compounds. \new{When the decision is made, the agent already knows whether the next lap will be run under clear track, VSC or SC, as rivals do.} In this section, we formally define the agent’s action space, state representation, reward function, and state transition dynamics.

\subsection{Action Space}
The agent’s action space consists of four possible actions: staying on track or changing tires to one of the three compounds available for the Grand Prix $g$, denoted by $\mathcal{C}^{(g,2024)}\subseteq \mathcal{C}$. 
\new{The coefficients of compounds rarely used at a circuit are obtained from the joint model described in Section~\ref{sec:cleartrack_base}, so all three compounds are available, and we write $\mathcal{C}^{r}=\mathcal{C}^{(g,2024)}$.} The action space is then defined as $\{0\}\cup\mathcal{C}^{r}$, where $0$ represents the action of staying on track (i.e., not performing a pit stop). \new{At the initial lap, the agent must select the starting tire compound, and action $0$ is penalized (Section~\ref{sec:reward}).} \new{Each driver starts the race with the dry tire sets available to them in the 2024 race, some of which were already used in practice or qualifying. Under the 2024 Sporting Regulations, a driver receives 13 dry sets per weekend (12 in sprint weekends) and returns some of them after each practice and qualifying session, which leaves about seven sets at the start of the race, including at least one hard and one medium set. Since the returned sets are not reported, we infer them from the sets used in each session (Appendix~\ref{app:tyre_sets}). A pit stop fits the least-worn set of the chosen compound, and the tire age of the new stint starts at the age of that set. The set taken off returns to the driver's sets with its age and can be fitted again, so a set of every compound is always available. Rivals fit their sets in the same way.}

\subsection{State Space}

The state space consists of features that describe relevant information about the agent and its opponents. These include: the current lap, the number of remaining laps, the tire compound of each car, the tire age of each car, \new{the age of the set of each compound that a pit stop would fit, }the time gap between the agent and each opponent, an indicator of whether each car has used at least two different tire compounds, the number of remaining laps under VSC or SC conditions, \new{the agent's action in the previous lap}, and the position of each car. A summary of the state variables is provided in Table~\ref{tab:estado_agente_individual}.

\begin{table}[htbp]
\centering
\begin{tabular}{>{\centering\arraybackslash}m{2cm}>{\centering\arraybackslash}m{3cm}>{\centering\arraybackslash}m{6.5cm}}
\toprule
Variable & Values & Description \\
\midrule
$S_{\text{Tire},d}$ & $\mathcal{C}^{r}$ & Tire compound in use of each driver $d\in\mathcal{D}$\\
$S_{\text{TireLife},d}$ & \new{$\mathbb{Z}_{\geq 0}$} & \new{Age in laps of the tire set in use of each driver $d\in\mathcal{D}$}\\
\new{$S_{\text{SetAge},c}$} & \new{$\mathbb{Z}_{\geq 0}$} & \new{Age of the set of compound $c\in\mathcal{C}^{r}$ that the agent would fit by stopping in the current lap}\\
$S_{\text{Lap}}$ & $\{1, \dots, N\}$ & Current lap number\\
$S_{\text{Laps Left}}$ & $\{1, \dots, N\}$ & Number of remaining laps in the race\\
$S_{\text{Change},d}$ & $\{0, 1\}$ & 1 if at least two different compounds have been used by driver $d\in\mathcal{D}$, 0 otherwise\\
$S_{\text{SC}}$ & $\{0, 1, 2, 3, 4\}$ & Number of remaining laps under SC\\
$S_{\text{VSC}}$ & $\{0, 1, 2\}$ & Number of remaining laps under VSC\\
\new{$S_{\text{PrevAction}}$} & $\{0\}\cup \mathcal{C}^{r}$ & \new{Action taken by the agent in the previous lap}\\
$S_{\text{Position},d}$ & $\{1, 2, 3, \dots, 20\}$ & Current position of car $d \in  \mathcal{D}$\\
$S_{\text{Gap},d}$ & \new{$\mathbb{R}$} & \new{Time difference in seconds between car $d \in \mathcal{D}$ and the agent, negative if $d$ is ahead}\\
\bottomrule
\end{tabular}
\caption{Description of state variables for the agent.}
\label{tab:estado_agente_individual}
\end{table}

\subsection{Rivals' Strategy}\label{sec:environment_rivals}

A key component of the environment is how the agent’s rivals make their decisions. On the one hand, rivals could act randomly by assigning probabilities to each possible action. On the other hand, they can follow more structured strategies, for instance by emulating decisions observed in real race data. The former approach may lead to unrealistic or suboptimal rival behavior, which can negatively affect the learning process. Therefore, we adopt the latter approach. Specifically, we use two predictive models to determine each rival’s action at every lap: the first model predicts whether the rival performs a pit stop or stays on track, and the second model selects the tire compound (from those available for the Grand Prix) if a pit stop is made. \new{During training, each rival may also retire from the race, with the probability and lap of retirement estimated from 2021--2023 races.} \new{When rivals' decisions are simulated, their starting compounds are drawn for each race from a common distribution $p\sim\text{Dirichlet}(\alpha,\alpha,\alpha)$ over the three available compounds, with $\alpha=0.437$ estimated by maximum likelihood from the starting compounds of 19 races held at the same circuits in 2021--2023, which reproduces the tendency of most drivers to start on the same compound.}

\subsubsection{Pit Stop Decision (LogitStop)}
We use a logistic regression model to estimate the probability that a driver performs a pit stop, as described by
\begin{equation}\label{eq:logit_rivals_stop_continue}
\begin{aligned}
\ln\left(\frac{q}{1-q}\right)=&\;\delta_{0}+\delta_{\text{SC}}X_{\text{SC}}+\delta_{\text{VSC}}X_{\text{VSC}}+\delta_{\text{Change}}X_{\text{Change}}\new{+\delta_{\text{Early}}X_{\text{Early}}+\delta_{\text{Late}}X_{\text{Late}}}\\
&+\sum_{c\in\mathcal{C}}\delta_{\text{TireLife},c}X_{\text{TireLife}}X_{c}+\delta_{\text{pos}}X_{\text{pos}}+\delta_{\text{GP}}X_{\text{GP}}.
\end{aligned}
\end{equation}
Here, $q$ denotes the probability of performing a pit stop. The variables $X_{\text{SC}}=\mathds{1}_{\{S_{\text{SC}}>0\}}$ and $X_{\text{VSC}}=\mathds{1}_{\{S_{\text{VSC}}>0\}}$ indicate whether the current lap is under SC or VSC conditions, respectively. \new{The variable $X_{\text{Change}}=S_{\text{Change}}$ indicates whether the driver has already used two compounds, while $X_{\text{Early}}=\mathds{1}_{\{S_{\text{Laps}}\leq 10\}}$ and $X_{\text{Late}}=\mathds{1}_{\{S_{\text{Laps}}\geq N-10\}}$ indicate the first and last ten laps of the race.} The variable $X_{\text{pos}}$ denotes the driver’s position, and $X_{\text{GP}}$ corresponds to a set of Grand Prix indicator variables. The variables $X_c$ and $X_{\text{TireLife}}$ are defined in Table~\ref{tab:variables_regression}. This model captures how race conditions and strategic context influence pit-stop decisions, and is referred to as \textbf{LogitStop}. \new{In the simulation, rivals do not stop in the first five laps and, under green flag conditions, neither with seven or fewer laps to go nor before completing eight laps on the current set; under SC or VSC, they may only stop on the first lap of the period. If a rival that has not yet used two compounds would fit the same compound again, a different compound is drawn at random.} Consistently with this rule, the model is estimated with one observation per driver and SC or VSC period, placed on the first lap of the period and equal to one if the driver stopped at any lap of it.

\subsubsection{Tire Compound Choice (CondLogitCompound)}

Once a pit-stop decision is made, the driver must choose the tire compound for the next stint. To model this decision, we use a conditional logistic regression calibrated using data from previous years. A conditional model is appropriate since the set of available compounds varies across year–Grand Prix combinations, allowing the model to account for the specific choice set at each pit stop. More specifically, the utility associated with selecting compound $c$ is specified as
\begin{equation}\label{eq:cond_logit_rivals}
U_c = \nu_{c} + \nu_{c,\text{LapsLeft}} X_{c} X_{\text{LapsLeft}} + \nu_{\text{NoChange}} X_{\text{NoChange}},
\end{equation}
where $c$ denotes the compound being evaluated, and $X_{\text{NoChange}}$ equals 1 if $c$ is the only compound that has been used so far. \new{$X_{\text{LapsLeft}}=N-S_{\text{Laps}}$ is the number of remaining laps, and compound C3 is the base alternative. Grand Prix indicators are not included, since they are common to all alternatives and cancel out in a conditional logit.} This conditional logistic model is denoted as \textbf{CondLogitCompound}.

\subsection{Reward Function}\label{sec:reward}

The reward in the Formula~1 simulation environment is defined as a combination of three components: (i) race time, (ii) race points, and (iii) penalties. The first component is given by the negative lap time (in seconds), computed using the regression models in Equations~\eqref{eq:reg} and~\eqref{eq:reg_yf}. This encourages the agent to minimize lap time and make decisions that improve race performance. The second component accounts for the points awarded based on the final race position. Points are assigned as 25, 18, 15, 12, 10, 8, 6, 4, 2, and 1 for the top ten positions, and 0 otherwise. To align the scale of this component with lap times, the points are multiplied by a factor of 300 before being incorporated into the reward. \new{In addition, $150\,(21-p)$ is added, where $p$ is the final position, so that the reward also distinguishes positions outside the top ten.} The third component is a penalty applied if the agent completes the race without using two different tire compounds. In this case (i.e., $S_{\text{Change},d}=0$), a penalty of \rev{1000} seconds is subtracted \new{and no points or position terms are added}. This ensures that the agent accounts for mandatory tire compound change rule. \new{Smaller penalties discourage stops that are rarely sensible in practice: 200 for stopping in two consecutive laps, 100 for stopping in laps 1 to 5, and 150 for stopping in the final lap. }\new{ Choosing action 0 at the initial lap, when a starting compound must be selected, assigns the medium compound and incurs a penalty of 1000.}

\rev{The inclusion of both race time and points reflects the trade-off between intermediate performance and final outcomes. While points determine the ultimate objective, \new{race time and the position term provide} a more informative signal in situations where positions result in identical points (e.g., positions 11 and 20). Reward normalization is applied \new{for PPO and Recurrent PPO} to improve training stability and convergence.}

\subsection{State Transition Function}

The state transition function in the Formula~1 simulation environment specifies how the state evolves as a result of the agent’s action $a \in \{0\}\cup\mathcal{C}^{r}$ for driver $d \in \mathcal{D}$. In particular, it defines the update rules for each relevant state variable, as detailed below.
\begin{align*}
S_{\text{Tire},d} &\leftarrow S_{\text{Tire},d}\cdot \mathds{1}_{\{a=0\}}+a\cdot \mathds{1}_{\{a\neq0\}} \\
S_{\text{TireLife},d} &\leftarrow (S_{\text{TireLife},d}+\mathds{1}_{\{S_{\text{SC}}+S_{\text{VSC}}=0\}})\cdot \mathds{1}_{\{a=0\}}+\new{W_{d,a}}\cdot \mathds{1}_{\{a\neq0\}} \\
S_{\text{Lap}} &\leftarrow S_{\text{Lap}}+1 \\
S_{\text{Laps Left}} &\leftarrow S_{\text{Laps Left}}-1 \\
S_{\text{Change},d} &\leftarrow \min\{S_{\text{Change},d}+\mathds{1}_{\{a\not\in\{0,S_{\text{Tire},d}\}\}},1\} \\
S_{\text{SC}} &\leftarrow 4\cdot\mathds{1}_{\{\text{SC begins}\}}+\max\{S_{\text{SC}}-1,0\} \\
S_{\text{VSC}} &\leftarrow 2\cdot\mathds{1}_{\{\text{VSC begins}\}}+\max\{S_{\text{VSC}}-1,0\} \\
\new{S_{\text{PrevAction}}} &\leftarrow a \\
S_{\text{Position},d} &\leftarrow \sum_{d'\in\mathcal{D}}\mathds{1}_{\{S_{\text{Gap},d}+y_{d}\le S_{\text{Gap},d'}+y_{d'}\}} \\
S_{\text{Gap},d} &\leftarrow S_{\text{Gap},d}+y_{d}-y_{\text{agent}}
\end{align*}
\new{where $W_{d,a}$ is the age of the least-worn set of compound $a$ among the sets of driver $d$, including the one just taken off. The set ages $S_{\text{SetAge},c}$ are updated accordingly.}
This set of rules fully specifies the evolution of the state variables as a function of the agent’s actions and race events, thereby defining the transition dynamics of the simulation environment.
\subsection{Stochastic Race Events}\label{sec:stochastic_race_events}
Besides the uncertainty incorporated in lap times and overtaking maneuvers, we also model the stochastic occurrence of VSC and SC events. Specifically, for each lap, we assume that a VSC or SC may be triggered with a certain probability. As noted in \cite{Aguad_Thraves}, these probabilities depend on the circuit characteristics. In particular, street circuits, which allow less margin for driver error, are more prone to yellow flag events, whereas more open circuits tend to exhibit lower probabilities. For each circuit, we estimate the probabilities of SC and VSC occurrences using empirical frequencies from previous seasons, see Appendix~\ref{app:prob_yf} for more details.

\subsection{Simulation Overview}
\begin{figure}[H]
\centering
\begin{tikzpicture}[
    every node/.style={font=\fontsize{7}{8}\sffamily\selectfont},
    model/.style={draw, rounded corners=3pt, minimum width=2.2cm, minimum height=0.65cm, align=center, fill=teal!20, draw=teal!50},
    simcomp/.style={draw, rounded corners=3pt, minimum width=2.5cm, minimum height=0.7cm, align=center, fill=purple!15, draw=purple!40},
    agent/.style={draw, rounded corners=3pt, minimum width=1.6cm, minimum height=0.65cm, align=center, fill=orange!20, draw=orange!40},
    arrow/.style={->, >=stealth, thin, gray!70},
    dashedborder/.style={draw, dashed, rounded corners=6pt, gray!60, inner sep=0.35cm}
]

\newcommand{\cylinder}[4]{%
    \begin{scope}[shift={(#1,#2)}]
        \fill[#4!15, draw=#4!40, line width=0.4pt] (-1.1,-0.38) rectangle (1.1,0.38);
        \fill[#4!25, draw=#4!40, line width=0.4pt] (0, 0.38) ellipse (1.1 and 0.22);
        \fill[#4!10, draw=#4!40, line width=0.4pt] (0,-0.38) ellipse (1.1 and 0.22);
        \node[font=\fontsize{7}{8}\sffamily\selectfont, align=center] at (0,0) {#3};
    \end{scope}
}

\cylinder{0}{-0.45}{Current GP \\ qualifying}{violet}
\node[inner sep=0, minimum width=2.2cm, minimum height=1.1cm] (qual) at (0,-0.45) {};

\cylinder{0}{-1.9}{Current GP \\ free practice}{violet}
\node[inner sep=0, minimum width=2.2cm, minimum height=1.1cm] (fp) at (0,-1.9) {};

\cylinder{0}{-3.35}{Previous \\ race data}{violet}
\node[inner sep=0, minimum width=2.2cm, minimum height=1.1cm] (prev) at (0,-3.35) {};

\node[model] (regltc)    at (3.8,  0)    {RegLTC};
\node[model] (regltyf)   at (3.8, -0.95) {RegLTYF};
\node[model] (logovert)  at (3.8, -1.9)  {LogitOvertake};
\node[model] (logstop)   at (3.8, -2.85) {LogitStop};
\node[model] (condlogit) at (3.8, -3.8)  {CondLogit\\Compound};

\node[simcomp] (laptime) at (7.6, -0.6)  {Lap time \\ simulation};
\node[simcomp] (reward)  at (7.6, -1.9)  {Reward function};
\node[simcomp] (rivals)  at (7.6, -3.2)  {Rivals' strategies};

\node[draw, dashed, rounded corners=6pt, black!75, inner sep=0.35cm,
      fit=(laptime)(reward)(rivals),
      label={[font=\fontsize{7}{8}\sffamily\selectfont, black!75]above:Simulation environment}] (simenv) {};
      
\node[agent] (agent) at (11.0, -1.9) {Agent};

\draw[arrow] (qual.east) -- (regltc.west);

\draw[arrow] (fp.east) -- (regltc.west);
\draw[arrow] (fp.east) -- (regltyf.west);
\draw[arrow] (fp.east) -- (logovert.west);

\draw[arrow] (prev.east) -- (regltc.west);
\draw[arrow] (prev.east) -- (regltyf.west);
\draw[arrow] (prev.east) -- (logovert.west);
\draw[arrow] (prev.east) -- (logstop.west);
\draw[arrow] (prev.east) -- (condlogit.west);

\draw[arrow] (regltc.east)    -- (laptime.west);
\draw[arrow] (regltyf.east)   -- (laptime.west);
\draw[arrow] (logovert.east)  -- (laptime.west);
\draw[arrow] (logstop.east)   -- (rivals.west);
\draw[arrow] (condlogit.east) -- (rivals.west);

\draw[arrow] (laptime.east) -- (agent.west);
\draw[arrow] (reward.east)  -- (agent.west);
\draw[arrow] (rivals.east)  -- (agent.west);

\node[anchor=north] at (1.0, -4.55) {%
    \tikz\fill[violet!15, draw=violet!40, line width=0.4pt]
    (0,0) rectangle (0.4,0.3);\quad Data sources};

\node[anchor=north] at (4.5, -4.55) {%
    \tikz\fill[teal!20, draw=teal!50, line width=0.4pt]
    (0,0) rectangle (0.4,0.3);\quad Predictive models};

\node[anchor=north] at (8.0, -4.55) {%
    \tikz\fill[purple!15, draw=purple!40, line width=0.4pt]
    (0,0) rectangle (0.4,0.3);\quad Sim. components};

\node[anchor=north] at (11.0, -4.55) {%
    \tikz\fill[orange!20, draw=orange!40, line width=0.4pt]
    (0,0) rectangle (0.4,0.3);\quad Agent};

\end{tikzpicture}
\caption{Overview of the data sources, predictive models, and simulation components used to construct the environment.}
\label{fig:summary}
\end{figure}
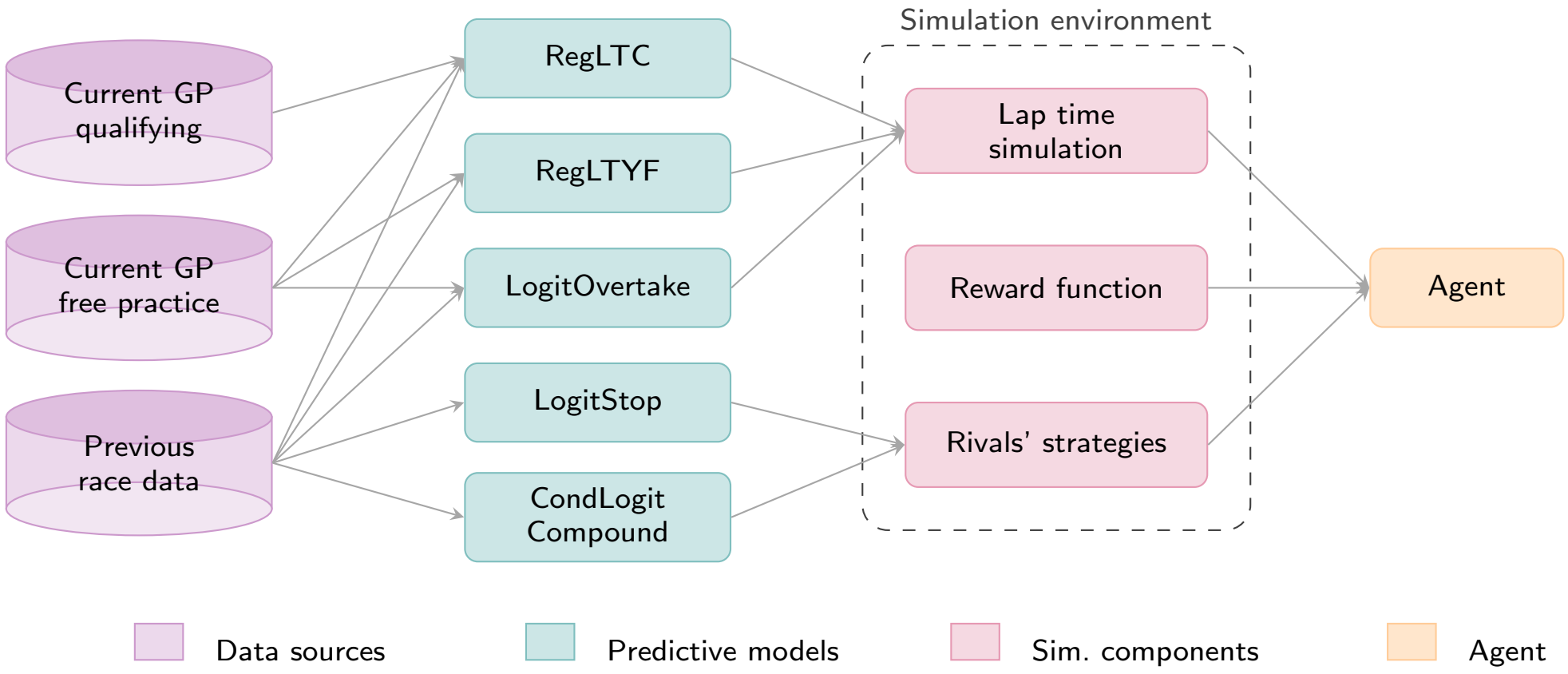

Figure~\ref{fig:summary} provides an overview of the complete simulation procedure. The predictive models are calibrated using historical data prior to the target Grand Prix, including past races as well as free-practice and qualifying sessions. These models, together with the reward function and rival behavior models, define the simulation environment used to train the agent. A detailed description of the lap-by-lap race evolution is provided in Algorithm~\ref{alg:simevolution} in Appendix~\ref{app:race_evolution}.

One agent is trained for each circuit, for a specific driver who starts from their actual 2024 grid position. Each training episode is a full race of $N$ laps in which lap times, overtaking maneuvers, rivals' decisions and retirements, and SC and VSC events are sampled as described above, and the agent chooses its starting compound. For each algorithm, several seeds are trained and one is selected on independent validation races. The hyperparameters, network architectures and selection rule for each DRL algorithm are reported in Appendix~\ref{app:architecture}.

\section{Exact Optimization Benchmark}\label{sec:QMIP}
This section describes an exact formulation for the pit-stop optimization problem in a deterministic single-driver setting. The formulation is based on quadratic mixed-integer programming (QMIP) and obtains optimal solutions when tire performance and pit stop time losses are known with certainty and no competition is present.

We adopt two simplifying assumptions consistent with the DRL setting. First, each tire compound $c$ exhibits a linear degradation rate, so that the lap time penalty grows linearly with tire usage. Under this assumption, the total lap time accumulated during a stint of $j$ laps \new{with a set of compound $c$ of age $w$ equals $\left(\beta_{c,1}+(w-1/2)\beta_{c,2}\right)j+\left(\beta_{c,2}/2\right)j^2$}. Second, pit stop time losses depend on track conditions: $\beta_{\text{PitIn}}+\beta_{\text{PitOut}}$ under clear track conditions\rev{, and $\beta_{\text{PitIn}}'+\beta_{\text{PitOut}}'$ under a VSC or a SC}. Although the driver cannot anticipate future yellow flag events, it is assumed to be aware of the current track status, denoted by $(S_{\text{VSC}},S_{\text{SC}})\in\{(0,0),(1,0),(0,1)\}$.

We define a \textit{stint} as a sequence of consecutive laps on the same tire compound. Consider that there are $n \leq N$ laps remaining, the current tire compound in use is $c^{(0)}$ \new{with age $w_{0}$}, and $\mathcal{C}^{(\mathrm{used})}$ denotes the set of compounds used in previous laps. The driver can complete at most $K$ additional stints. \new{Let $w_{1c}\leq w_{2c}\leq\dots$ be the ages of the sets of compound $c$ left, of which at most $K$ are considered, and let $\mathcal{J}_c$ index them. The decision variables are $x_{0}\geq 0$, the laps remaining in the current stint, and $x_{jc}\geq 0$, the laps run on the $j$-th least-worn set of compound $c$, with binary indicator $y_{jc}$.} Following \cite{QMIP}, the problem is formulated as:
\begin{subequations}\label{eq:QMIP}
\begin{alignat}{3}
\min \quad & \new{\left(\beta_{c^{(0)},1}+\left(w_{0}-\frac{1}{2}\right)\beta_{c^{(0)},2}\right)x_{0}+\frac{\beta_{c^{(0)},2}}{2}x^2_{0}} && \notag\\
& \mathrlap{\new{+\sum_{c \in \mathcal{C}^{r}}\sum_{j\in\mathcal{J}_c}\left[\left(\beta_{c,1}+\left(w_{jc}-\frac{1}{2}\right)\beta_{c,2}\right)x_{jc}+\frac{\beta_{c,2}}{2}x^2_{jc}+a\,y_{jc}\right]}} \label{eq:QMIP_obj}\\
\text{s.t.} \quad & \new{x_{0}+\sum_{c \in \mathcal{C}^{r}}\sum_{j\in\mathcal{J}_c} x_{jc} = n} && \label{eq:QMIP_a}\\
           & \new{x_{jc} \geq y_{jc}} 
             && \new{\forall c \in \mathcal{C}^{r},\;\forall j \in \mathcal{J}_c} \label{eq:QMIP_c}\\
           & \new{x_{jc} \leq n \cdot y_{jc}} 
             && \new{\forall c \in \mathcal{C}^{r},\;\forall j \in \mathcal{J}_c} \label{eq:QMIP_b}\\
           & \new{\sum_{c \in \mathcal{C}^{r}\setminus\mathcal{C}^{(\mathrm{used})}}\sum_{j\in\mathcal{J}_c} y_{jc} \geq 1} 
             && \quad\text{if } \left|\mathcal{C}^{(\mathrm{used})}\right| = 1 \label{eq:QMIP_e}\\
           & \new{\sum_{c \in \mathcal{C}^{r}}\sum_{j\in\mathcal{J}_c} y_{jc}\leq K} && \label{eq:QMIP_f}\\
           & \new{y_{jc}\geq y_{j+1,c}} 
             && \new{\forall c \in \mathcal{C}^{r},\;\forall j,j+1 \in \mathcal{J}_c} \label{eq:QMIP_h}\\
           & \new{x_{0},x_{jc} \in \mathbb{Z}_{\geq 0}} 
             && \new{\forall c \in \mathcal{C}^{r},\;\forall j \in \mathcal{J}_c} \label{eq:QMIP_n}\\
           & \new{y_{jc}\in \{0,1\}} 
             && \new{\forall c \in \mathcal{C}^{r},\;\forall j \in \mathcal{J}_c} \label{eq:QMIP_o}
\end{alignat}
\end{subequations}
where \rev{$a_{1}=(\beta_{\mathrm{PitIn}}+\beta_{\mathrm{PitOut}})\mathds{1}_{\{S_{\mathrm{SC}}+S_{\mathrm{VSC}}=0\}}+(\beta_{\mathrm{PitIn}}'+\beta_{\mathrm{PitOut}}')\mathds{1}_{\{S_{\mathrm{SC}}+S_{\mathrm{VSC}}>0\}}$} and \new{$a=\beta_{\mathrm{PitIn}}+\beta_{\mathrm{PitOut}}$} denote the time loss incurred during a pit stop under current and normal track conditions, respectively. \new{The reduced loss $a_{1}$ replaces $a$ in one stop if the driver stops in the current lap ($x_{0}=0$), which is modeled with an auxiliary binary variable; under SC, this is only possible in its first lap.} The objective function~\eqref{eq:QMIP_obj} minimizes the total remaining race time. Constraint~\eqref{eq:QMIP_a} ensures that the number of laps assigned across all stints equals the remaining laps. Constraints~\eqref{eq:QMIP_c} and~\eqref{eq:QMIP_b} enforce the logical relationship between the integer and binary decision variables. Constraint~\eqref{eq:QMIP_e} guarantees that at least two different tire compounds are used throughout the race, as required by regulations. \new{Since the order of the stints after the current one does not affect the objective, they are indexed by tire set rather than by position, which removes equivalent permutations: each set is fitted at most once, at most $K$ sets are fitted (Constraint~\eqref{eq:QMIP_f}), and the sets of a compound are fitted from the least worn (Constraint~\eqref{eq:QMIP_h}). Two further constraints remove equivalent solutions for the current stint: at the start of the race at least one lap is run on the starting compound and, without a yellow flag, the least-worn set of the current compound is not run for more laps than the current set will have completed.} 

The formulation can be readily adapted to optimize pit-stop strategies from the start of the race by setting the initial tire compound as an additional decision variable\new{, which starts on the least-worn set of that compound}. In deterministic instances without yellow flag events, the QMIP~\eqref{eq:QMIP} is solved once at the start of the race to obtain the full strategy. In stochastic instances where SC or VSC events may occur, the QMIP~\eqref{eq:QMIP} is solved at each lap, allowing the strategy to be updated in response to the current track status. \new{The QMIP does not refit a set taken off during the race. In the deterministic instances of Section~\ref{sec:results}, an exact dynamic program over stints gives the same optimal race time with and without refitting sets.}

\section{Validation of Learning-Based Pit-Stop Policies}\label{sec:results}

\rev{Before using the learned policies as decision-support tools, we first evaluate their reliability in settings of increasing complexity. This section is organized into three parts. First, we present the calibration of the simulation environment, including the estimated models for lap times, driver coefficients, and overtaking behavior. Second, we benchmark the DRL methods against the optimal solution obtained via the QMIP formulation in simplified settings without uncertainty or competition. Finally, we evaluate the performance of the DRL methods and the QMIP in realistic race scenarios with a full grid of competitors, considering both fully simulated instances and instances based on real race realizations.}

\subsection{Calibration of the Simulation Environment}

\subsubsection{Clear-Track Lap-Time Model}

\begin{figure}[H]
\centering
\includegraphics[width=\linewidth]{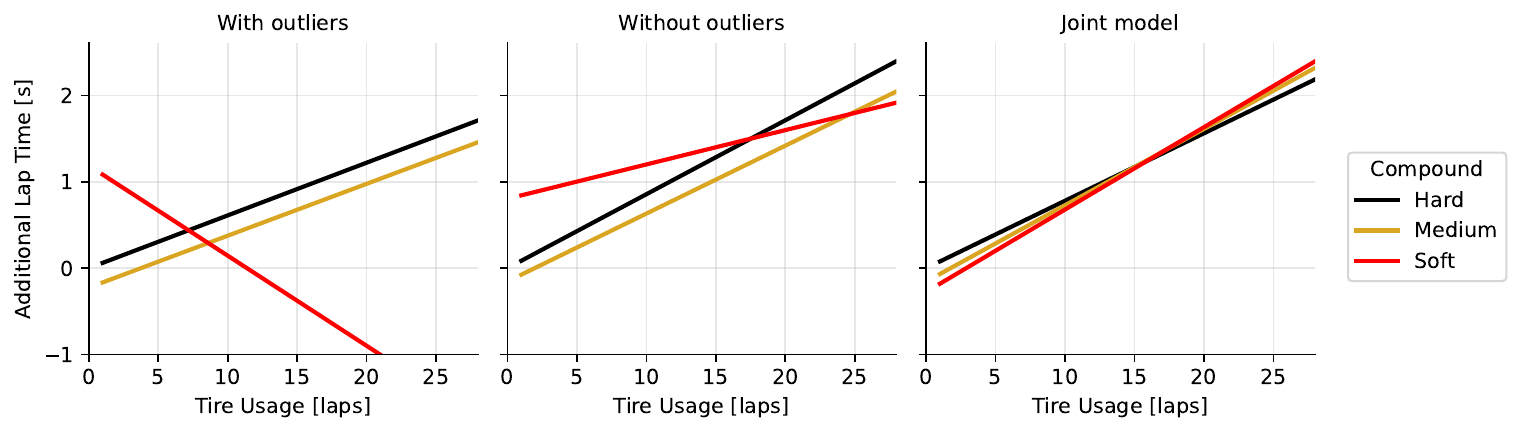}
\caption{Tire performance at the 2024 United States GP. Left: Equation~\eqref{eq:reg} estimated with all laps. Center: Equation~\eqref{eq:reg} estimated after removing outliers with Algorithm~\ref{alg:IQR}. Right: coefficients of the joint model (Appendix~\ref{app:regTC_tirecoef_matrix_completion}). Additional lap time is relative to a new hard tire.}
\label{fig:tire_compounds}
\end{figure}

\new{The regression model in Equation~\eqref{eq:reg} is estimated separately for each circuit, following the outlier removal procedure described in Algorithm~\ref{alg:IQR}, which removes between 0.7\% and 3.0\% of the laps in the circuits considered. Figure~\ref{fig:tire_compounds} shows the resulting tire compound coefficients for the 2024 United States Grand Prix. With all laps, the soft compound becomes faster as it ages. Removing outliers restores a positive degradation, but the soft compound remains slower than the hard one when new and degrades more slowly than the other two, since it was used in only five stints in 2021--2023. The joint model of Appendix~\ref{app:regTC_tirecoef_matrix_completion} recovers the expected pattern, in which softer compounds are faster when new but degrade faster.}

\new{Table~\ref{tab:regression_cleanlap} in Appendix~\ref{app:table_reg_clean} reports the full set of estimated coefficients before and after removing outliers for the 2024 United States Grand Prix, and Figure~\ref{fig:laptimes_outliers} in Appendix~\ref{app:residuals} shows the removed lap times.}

\new{The estimated coefficients of the yellow flag lap time model (RegLTYF) are reported in Appendix~\ref{app:regLTYF}. A lap under VSC is about $34\%$ slower than a clear-track lap, and a pit stop under yellow flag conditions costs between $59\%$ and $69\%$ of the clear-track pit-lane loss across the eight circuits.}

\subsubsection{Driver Coefficient Correction}
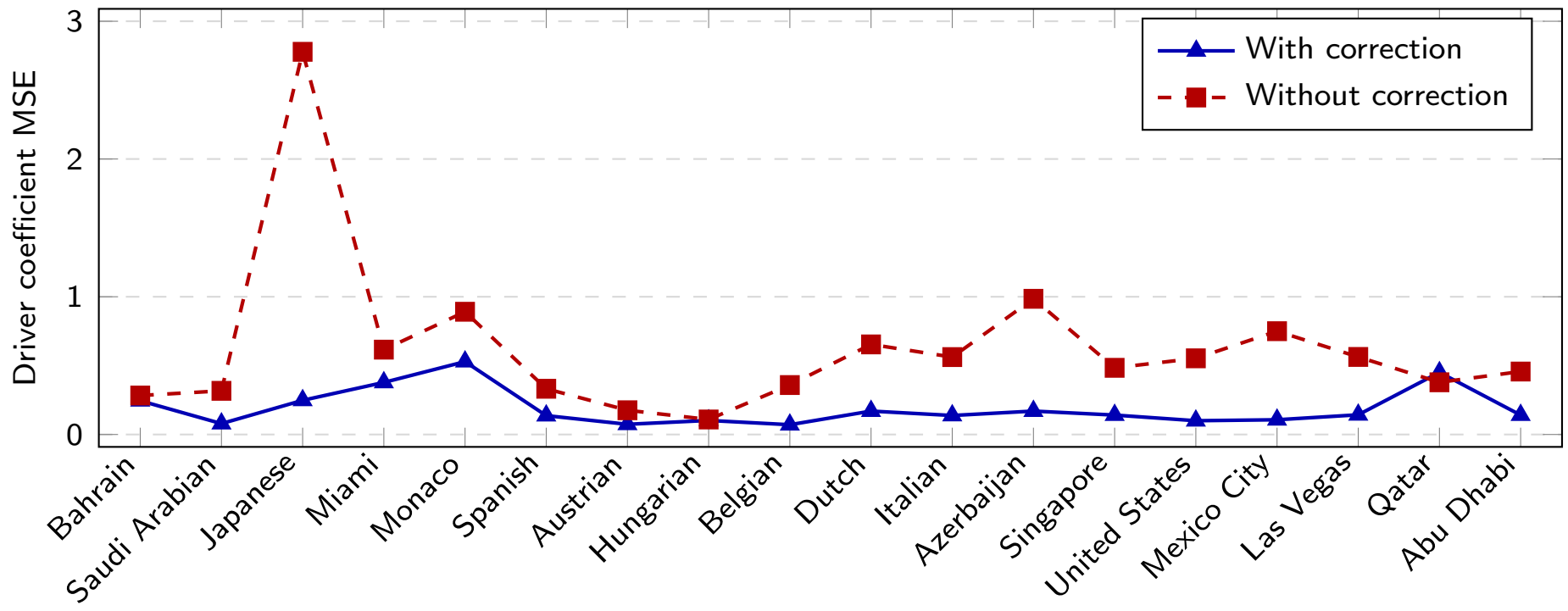
\begin{figure}[H]
\centering
\begin{tikzpicture}
\begin{axis}[
    width=13cm,
    height=5cm,
    ylabel={Driver coefficient MSE},
    xtick=data,
    xticklabels={
        Bahrain,
        Saudi Arabian,
        Japanese,
        Miami,
        Monaco,
        Spanish,
        Austrian,
        Hungarian,
        Belgian,
        Dutch,
        Italian,
        Azerbaijan,
        Singapore,
        United States,
        Mexico City,
        Las Vegas,
        Qatar,
        Abu Dhabi
    },
    xticklabel style={rotate=45, anchor=east, font=\footnotesize\sffamily},
    ymin=0,
    ymax=3.0,
    ymajorgrids=true,
    grid style={dashed, gray!30},
    legend style={
        at={(0.98,0.98)},
        anchor=north east,
        font=\footnotesize\sffamily,
        cells={anchor=west}
    },
    tick label style={font=\footnotesize\sffamily},
    yticklabel={\pgfmathprintnumber[assume math mode=true]{\tick}}, 
    label style={font=\footnotesize\sffamily},
    enlargelimits=0.03,
]
\addplot[
    color=blue!70!black,
    mark=triangle*,
    mark size=1.8pt,
    thick,
] coordinates {
    (1, 0.245)
    (2, 0.079)
    (3, 0.248)
    (4, 0.378)
    (5, 0.529)
    (6, 0.137)
    (7, 0.074)
    (8, 0.102)
    (9, 0.071)
    (10, 0.170)
    (11, 0.138)
    (12, 0.170)
    (13, 0.141)
    (14, 0.100)
    (15, 0.107)
    (16, 0.143)
    (17, 0.451)
    (18, 0.141)
};
\addlegendentry{With correction}
\addplot[
    color=red!70!black,
    mark=square*,
    mark size=1.8pt,
    thick,
    dashed,
    mark options={solid} 
] coordinates {
    (1, 0.283)
    (2, 0.317)
    (3, 2.778)
    (4, 0.616)
    (5, 0.891)
    (6, 0.332)
    (7, 0.176)
    (8, 0.109)
    (9, 0.359)
    (10, 0.654)
    (11, 0.562)
    (12, 0.985)
    (13, 0.484)
    (14, 0.553)
    (15, 0.750)
    (16, 0.563)
    (17, 0.379)
    (18, 0.457)
};
\addlegendentry{Without correction}
\end{axis}
\end{tikzpicture}
\caption{MSE between predicted and ground-truth driver coefficients for each 2024 Grand Prix, where the ground truth is estimated with Equation~\eqref{eq:reg} using only that race. \textit{Without correction}: Equation~\eqref{eq:reg} on the 2021--2023 editions of the same circuit; \textit{With correction}: Equation~\eqref{eq:regcoef}, which adds the driver's three most recent 2024 races plus practice and qualifying sessions of the same weekend. Circuits without dry editions in 2024 or before are omitted, as is the Australian GP, where the reference driver retired on the fourth lap.}
\label{fig:mse_comparison}
\end{figure}

\new{Figure~\ref{fig:mse_comparison} shows the prediction accuracy of the driver coefficients with and without the correction of Equation~\eqref{eq:regcoef}. Incorporating recent races and the practice and qualifying sessions of the same weekend reduces the error in all circuits except the Qatar GP. An example of the estimated coefficients of Equation~\eqref{eq:regcoef} is shown in Appendix~\ref{app:updatebelgium}.}

\subsubsection{Overtaking Model}
The logistic regression in Equation~\eqref{eq:logitovertake} is estimated using data from previous years. Table~\ref{tab:overtake_logit} in Appendix~\ref{app:overtake} reports the estimated coefficients for each circuit. In all cases, the intercept is negative, indicating that closing the gap to the car ahead does not necessarily lead to an overtake. \new{The coefficient on the predicted time gap is negative and significant at all circuits, indicating that smaller gaps are associated with higher overtaking probabilities. The coefficients vary across circuits, reflecting differences in overtaking opportunities: overtaking is more frequent at permanent circuits such as Belgium or Bahrain, and less frequent at street circuits such as Singapore or Saudi Arabia.}

\subsubsection{Rivals' Strategy Models}
The logistic regression for the pit-stop decision (LogitStop) and the conditional logistic regression for tire compound selection (CondLogitCompound) are estimated using historical race data as described in Section~\ref{sec:environment_rivals}. For space considerations, the estimated coefficients are reported in Appendix~\ref{app:rivals_strategy}.

\subsection{Optimality Gap of DRL Methods in Simplified Settings}
To measure how close each DRL method comes to the optimum, we evaluate the agent on instances simple enough so that the optimal pit-stop strategy can be computed. More specifically, we consider the agent for a specific GP under no competition in a deterministic setting where tires' performance and pit-stop additional lap time are obtained from the regression in Equation~\eqref{eq:reg}. \new{The agent has the tire sets of the driver considered at the same circuit in Section~\ref{sec:agent_competition}.} Then, the estimated coefficients ($\{\beta_{c,1},\beta_{c,2}\}_{c\in\mathcal{C}}$, $\beta_{\text{PitIn}}$ and  $\beta_{\text{PitOut}}$) are used as inputs to the QMIP in Equation~\eqref{eq:QMIP} of Section~\ref{sec:QMIP} to obtain the optimal solution. \new{Each DRL method is trained with five seeds for each circuit, and we report the seed with the smallest gap among those that satisfy the two-compound rule.}

\begin{figure}[H]
\centering
\includegraphics[width=1\linewidth]{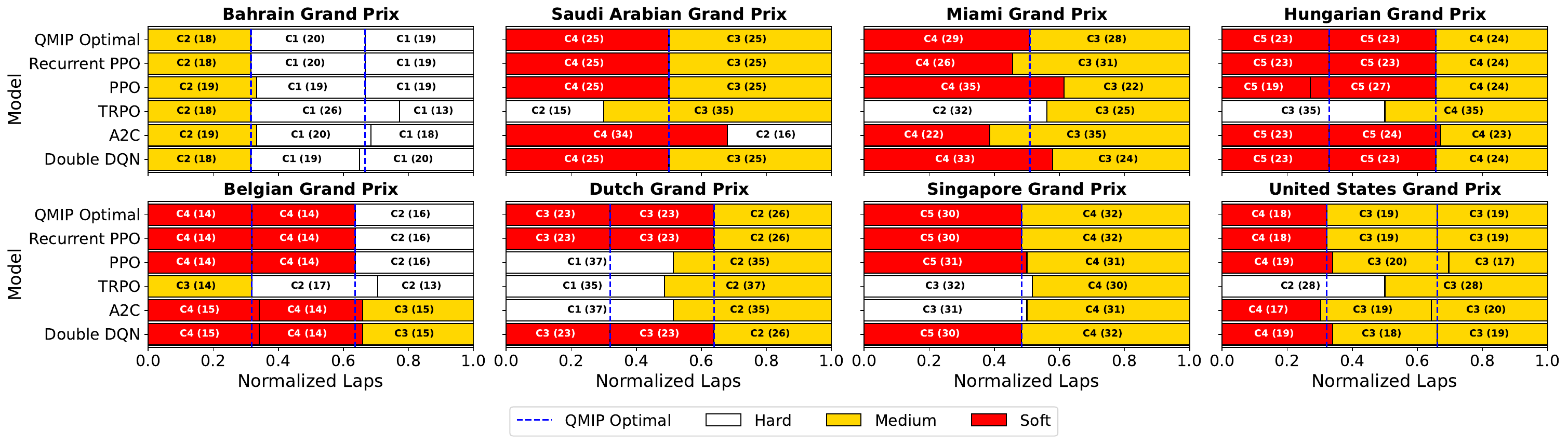}
\caption{Pit-stop strategies for five DRL models and the optimal solution, in a deterministic setting with no competition.}
\label{fig:stints_alone}
\end{figure}

Figure~\ref{fig:stints_alone} shows, for each GP, the optimal strategy, obtained with the QMIP, and the agent's strategy with all DRL methods. In general, the pit-stop solutions obtained with all methods, except with TRPO, are very similar to the optimal solution. With the selected seed, Recurrent PPO attains the optimal race time in seven of the eight circuits and Double DQN in five. Since different strategies can have the same race time, we also report the race time suboptimality gap in Table~\ref{tab:opt_gap} in Appendix~\ref{app:opt_gap_simple_instance}. Averaged across circuits, the gap of the selected seed is below one second for Recurrent PPO, Double DQN and PPO, about one and a half seconds for A2C, and larger for TRPO. With the median over seeds, the gaps increase for all methods: they remain below one second for Double DQN and Recurrent PPO, rise to between two and three seconds for A2C and PPO, and exceed five seconds for TRPO.

\subsection{Benchmarking DRL and QMIP under Competition}\label{sec:agent_competition}

This section analyzes the performance of the DRL methods and the QMIP on \new{eight} circuits of the 2024 season, considering that the agent competes in a full grid of twenty drivers. For each GP, we consider a specific driver \new{starting between sixth and thirteenth on the grid}, so that there is room for both improvement and deterioration in the final result. Once the agent has been trained for a particular GP in a simulator calibrated with data from previous races, using the seed selected on validation races (Figure~\ref{fig:seed_spread} in Appendix~\ref{app:architecture} shows all seeds), we evaluate its performance over 500 simulated races under two settings: one in which rivals' decisions and yellow flag events are randomly generated (as described in Sections~\ref{sec:environment_rivals} and \ref{sec:stochastic_race_events}, respectively), and one in which they are taken from the actual 2024 race\new{, including the age of the set fitted in each stint}. \new{In both settings, rivals that retired in the actual race retire on the same lap (Table~\ref{tab:DNFs_2024} in Appendix~\ref{app:DNF}). The QMIP uses the same tire coefficients as the DRL methods, estimated with 2021--2023 data, and is re-optimized at each lap given the current track status. In both settings, the clear-track lap times of the evaluation races are simulated with tire coefficients re-estimated with the 2021--2024 races, so the agents and the QMIP are evaluated in a simulator that differs from the one used to train the agents.}

\subsubsection{Fully Simulated Race Scenarios}\label{sec:competition_fully}

\begin{figure}[H]
\centering
\includegraphics[width=1\linewidth]{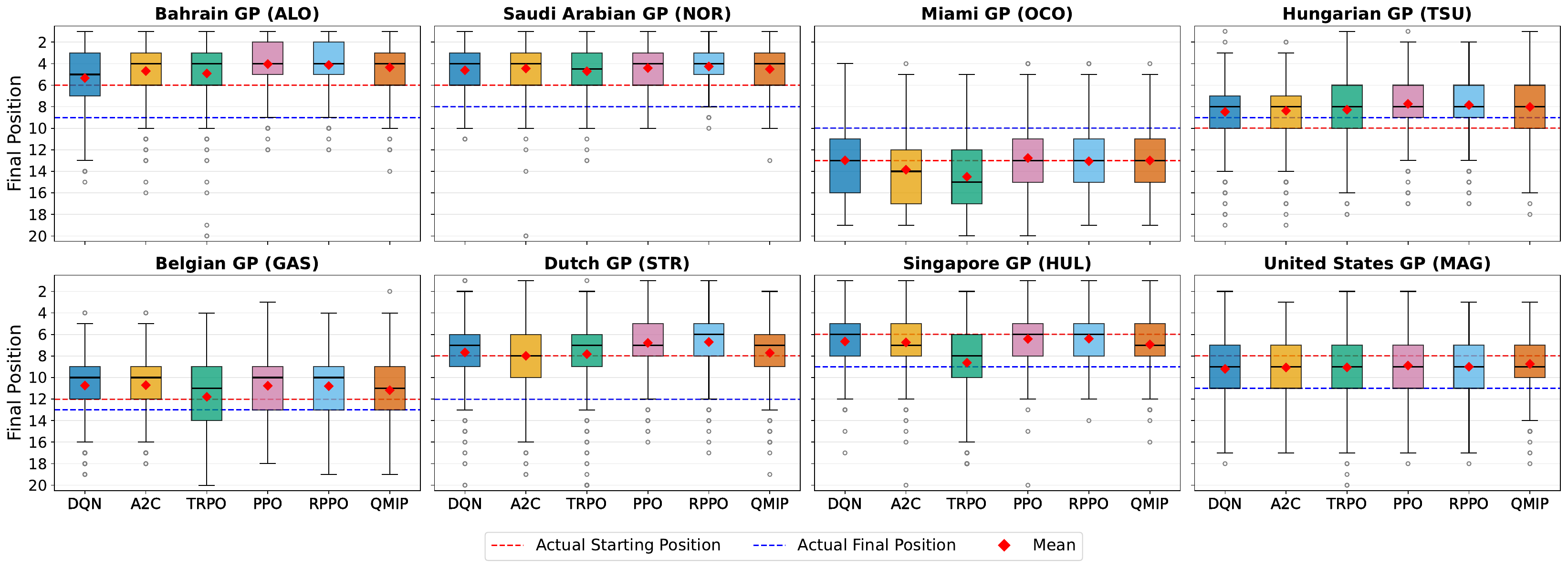}
\caption{Distribution of final positions over $500$ fully simulated races per circuit, 2024 season.}
\label{fig:stats_drivers_stochastic_2024}
\end{figure}

\new{Figure~\ref{fig:stats_drivers_stochastic_2024} shows the distribution of final positions across 500 simulations for each solution method and circuit in the 2024 season. In terms of the median, PPO and Recurrent PPO are never worse than the QMIP and are better in two and three circuits, respectively. The dispersion of final positions is similar across methods. Mean final positions with 95\% confidence intervals are reported in the top part of Table~\ref{tab:ci_2024} in Appendix~\ref{app:ci_2024}. Averaged over the eight circuits, PPO and Recurrent PPO attain the best mean final positions, with no appreciable difference between them, improving on the QMIP by $0.33$ and $0.28$ positions, respectively, followed by the QMIP, then Double DQN and A2C, and finally TRPO.}

\subsubsection{Real Race Realizations}\label{sec:competition_racescenario}

\begin{figure}[H]
\centering
\includegraphics[width=1\linewidth]{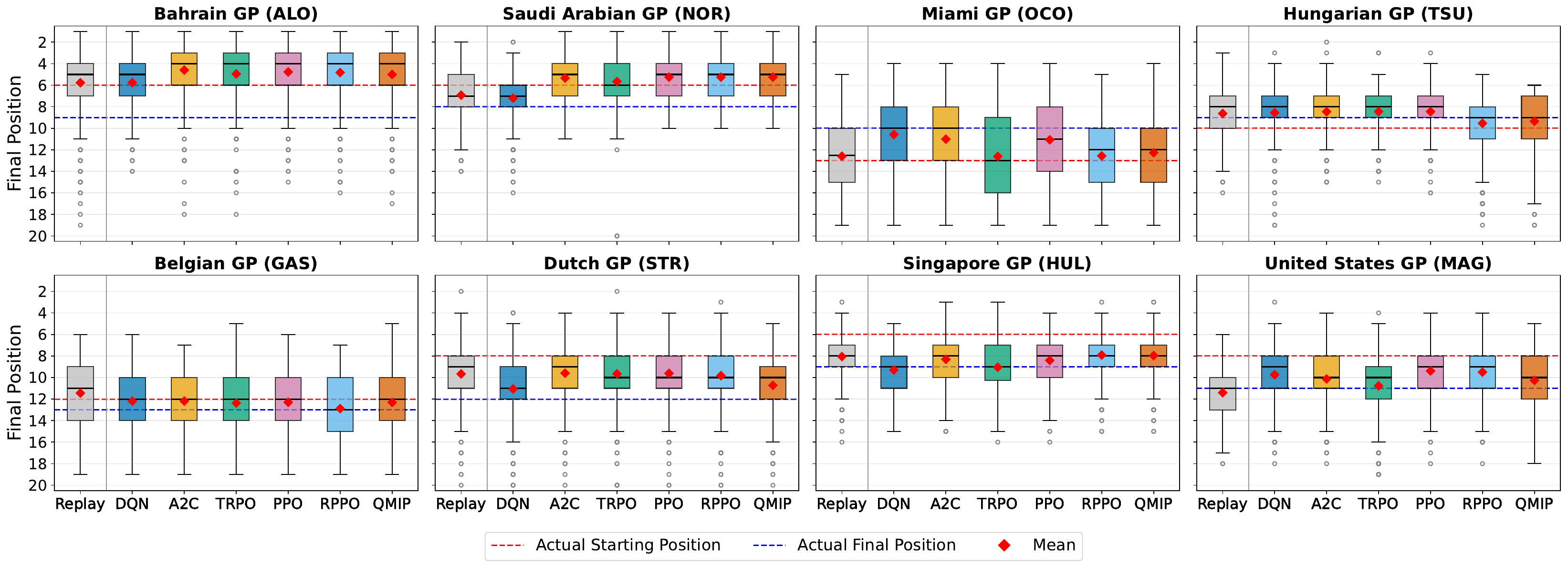}
\caption{Distribution of final positions over $500$ simulations per circuit with rivals' decisions and yellow flags fixed to those of the actual 2024 race. The leftmost box, ``Replay'', reproduces the driver's actual 2024 starting compound and pit stops.}
\label{fig:boxplot_drivers_2024}
\end{figure}

\new{Figure~\ref{fig:boxplot_drivers_2024} shows the distribution of final positions for each solution method under real race scenarios for the 2024 season. Compared to the fully simulated setting, the distributions are somewhat narrower on average, although not in every circuit, given that rivals' decisions and yellow flag events are fixed to those observed in the actual race. From the bottom part of Table~\ref{tab:ci_2024} in Appendix~\ref{app:ci_2024}, we observe that PPO and A2C obtain the best average final positions, $0.49$ and $0.44$ positions ahead of the QMIP, followed by Recurrent PPO and the QMIP, and finally TRPO and Double DQN. As a reference, the box labeled ``Replay'' in Figure~\ref{fig:boxplot_drivers_2024} reproduces the driver's actual 2024 starting compound and pit stops in the same simulated races. Averaged over the eight circuits, every method except Double DQN improves on the final position of this replay, by $0.66$ positions for PPO and $0.61$ for A2C. The replay, however, outperforms every method in the Belgian GP. The replay also provides a check of the simulator: the driver's actual final position (dashed line) lies within the box of the replayed positions in six of the eight races, and within their central $90\%$ interval in all of them. On average, however, the replay places the driver $0.94$ positions ahead of the actual final position, which is why the methods are compared with the replay rather than with the actual final position.}

\new{Across the eight circuits, PPO and Recurrent PPO are the only methods that improve on the QMIP in both settings, PPO by a larger margin. The differences are statistically significant (Table~\ref{tab:ci_2024}). These margins are per race: half a position in each of the 24 races of a season adds up to twelve positions. In the 2024 Constructors' Championship, Alpine finished sixth, seven points ahead of Haas \citep{f1_2024_standings}, and in 2023 sixth place was worth an estimated 8 million dollars more in prize money than seventh \citep{hardy2023prize}. We therefore adopt PPO for the look-ahead analysis, whose rollouts sample yellow flag events, the same uncertainty under which the agent was trained.}

\section{AI-Assisted Look-Ahead Decision Support}\label{sec:decision_support}

While the previous sections evaluate decision policies in terms of their average performance over many races, in practice a race strategist must make a single, irreversible choice on each lap, with limited time for deliberation. In this setting, a point recommendation from the DRL agent, such as ``pit for medium tires now'', offers little visibility into the underlying trade-offs, into how confident the policy is in that action, or into what would have happened under alternative choices. For example, if the race engineer believes that the best action is to ``continue on track'' rather than stopping for medium tires, it would be useful to visualize and compare the race outcomes implied by the DRL recommendation and by the engineer's own preferred action. We therefore propose a visualization tool that turns the trained DRL policy into a forward look-ahead simulator over alternative actions, providing the strategist with an interpretable view of the distribution of possible outcomes under each available action.

For a fixed race, lap $n$, and driver, the procedure is as follows. For each candidate action $a \in \{\text{Hard}, \text{Medium}, \text{Soft}, \text{No pit}\}$ available at lap $n$, we run $M$ Monte Carlo simulations in which the action at lap $n$ is forced to $a$ and the actions in laps $n+1, \dots, N$ \new{are those that the trained} \rev{PPO} \new{policy $\pi$ would take, that is, its most likely action, excluding stops for the same compound before the two-compound rule is met}. Each rollout draws independent realizations of the stochastic race events described in Section~\ref{sec:stochastic_race_events}\rev{ (yellow flag occurrences, lap time noise, and overtakes),} as well as of the rivals' pit stop decisions through the predictive model of Section~\ref{sec:environment_rivals}. The final position and championship points of each rollout give an empirical distribution of outcomes for each action at lap $n$. \new{The state $s_n$ is obtained by replaying the actual 2024 race up to lap $n$, with rivals' decisions and yellow flags fixed to those observed and the driver following $\pi$. A yellow flag period already active at lap $n$ keeps its actual duration, and later events are sampled. Since the policy is followed deterministically, the spread of outcomes reflects only race uncertainty and not randomness in the policy's own actions.} The procedure is summarized in Algorithm~\ref{alg:counterfactual_rollouts}. It is a one-step rollout with the PPO policy as the base policy \citep{bertsekas1997rollout}.

\begin{algorithm}[H]
\caption{Look-ahead rollouts for decision support visualization}\label{alg:counterfactual_rollouts}
\begin{algorithmic}[1]
\Require Trained policy $\pi$, current race state $s_n$, candidate actions $\{0\}\cup\mathcal{C}^{r}$, number of simulations $M$\new{, set of feasible actions $\mathcal{A}(s)$}
\For{each action $a \in \{0\}\cup\mathcal{C}^{r}$}
    \For{$m = 1, \dots, M$}
        \State Set $a_n \gets a$
        \State Sample stochastic race events and rivals' actions; update state $s_{n+1}$
        \For{$n' = n+1, \dots, N$}
            \State $a_{n'} \gets \arg\max_{a' \in \mathcal{A}(s_{n'})} \pi(a' \mid s_{n'})$
            \State Sample stochastic race events and rivals' actions
            \State Update state $s_{n'+1}$
        \EndFor
        \State \rev{Record final outcome $y_m^{(a)}$ (final position and championship points)}
    \EndFor
\EndFor
\State \Return Empirical distributions $\{y_m^{(a)}\}_{m=1}^{M}$ for each $a \in \{0\}\cup\mathcal{C}^{r}$
\end{algorithmic}
\end{algorithm}

\begin{figure}[H]
\centering
\includegraphics[width=0.95\linewidth]{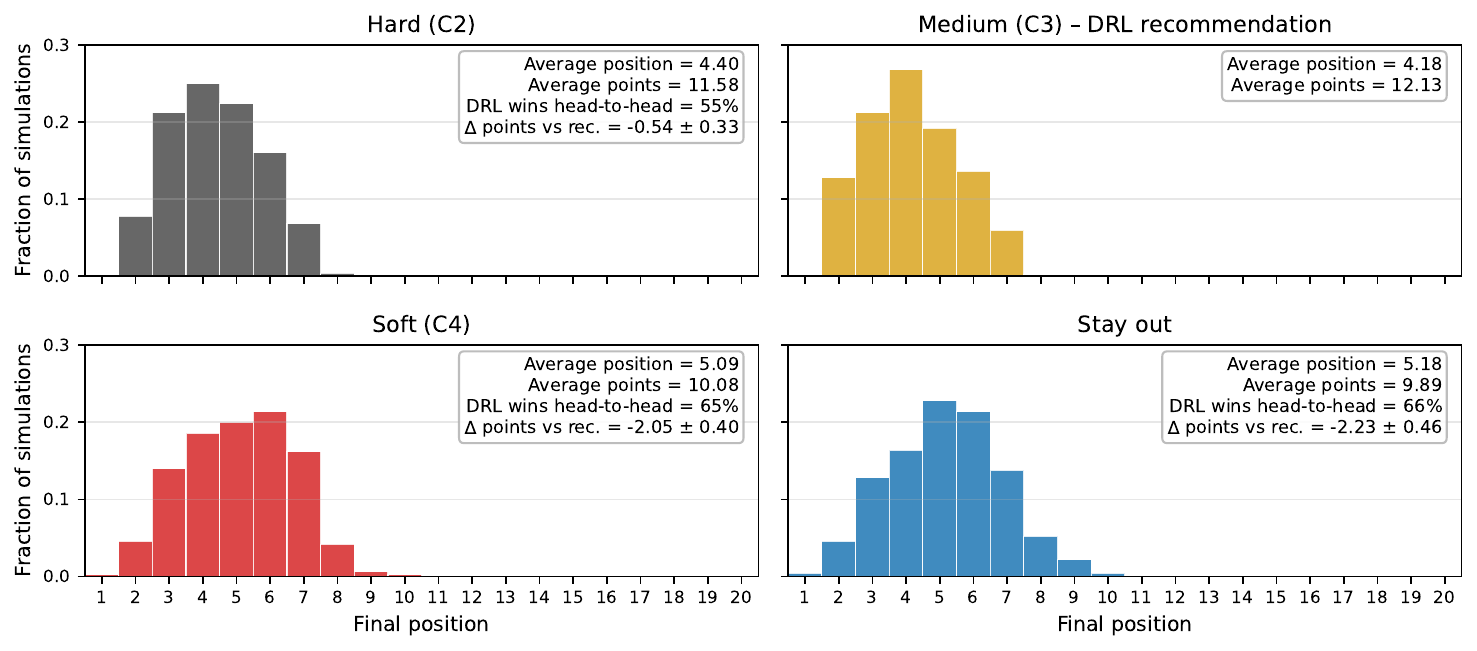}
\caption{Final position over $M = 500$ rollouts of each candidate action at lap $7$ of the 2024 Saudi Arabian Grand Prix, with the PPO policy driving in place of Lando Norris. Head-to-head wins and $\Delta$ points compare paired rollouts with the recommended action ($\pm 2$ standard errors, SE).}
\label{fig:decision_support}
\end{figure}

\new{Figure~\ref{fig:decision_support} illustrates the resulting visualization for a representative decision point: lap $n = 7$ of the 2024 Saudi Arabian Grand Prix ($N = 50$ laps), in a race driven by the PPO policy from Lando Norris's starting position. The policy started on soft tires and is running in position $6$, with seven laps of use on them, on the lap in which a safety car is deployed. The DRL policy recommends switching to medium tires, which attains the best average final position ($4.18$) and the highest average points ($12.13$)\footnote{Only drivers finishing in the top ten positions receive championship points.}. Fitting hard tires comes closest, with $0.54 \pm 0.33$ fewer points, so the strategist may prefer it at a small expected cost on the basis of information not encoded in the simulator, such as weather forecasts, mechanical concerns, or tactical instructions from the team principal. Staying on track loses $2.23 \pm 0.46$ points. A stop under the safety car costs less time than under green flag conditions and is only possible on the first lap of the period (Appendix~\ref{app:race_evolution}), so if the driver stays out, the policy itself postpones the stop to the second half of the race in about $70\%$ of the rollouts. Fitting soft tires again costs about as much ($2.05 \pm 0.40$ fewer points), because a second soft stint would still require a later stop on the medium or hard compound. Hard and medium tires complete the race with the stop at lap $7$, whereas fitting softs again or staying on track lead to a further stop (see Figure~\ref{fig:decision_stints} in Appendix~\ref{app:decision_stints}). The strategist can thus tell near-equivalent options from clearly dominated ones. The actual race differs: Norris started on medium tires and stayed out under the safety car, and the policy would have recommended a stop for soft tires at that lap (Figure~\ref{fig:realhist} in Appendix~\ref{app:realhist}).}

\begin{figure}[H]
\centering
\includegraphics[width=0.95\linewidth]{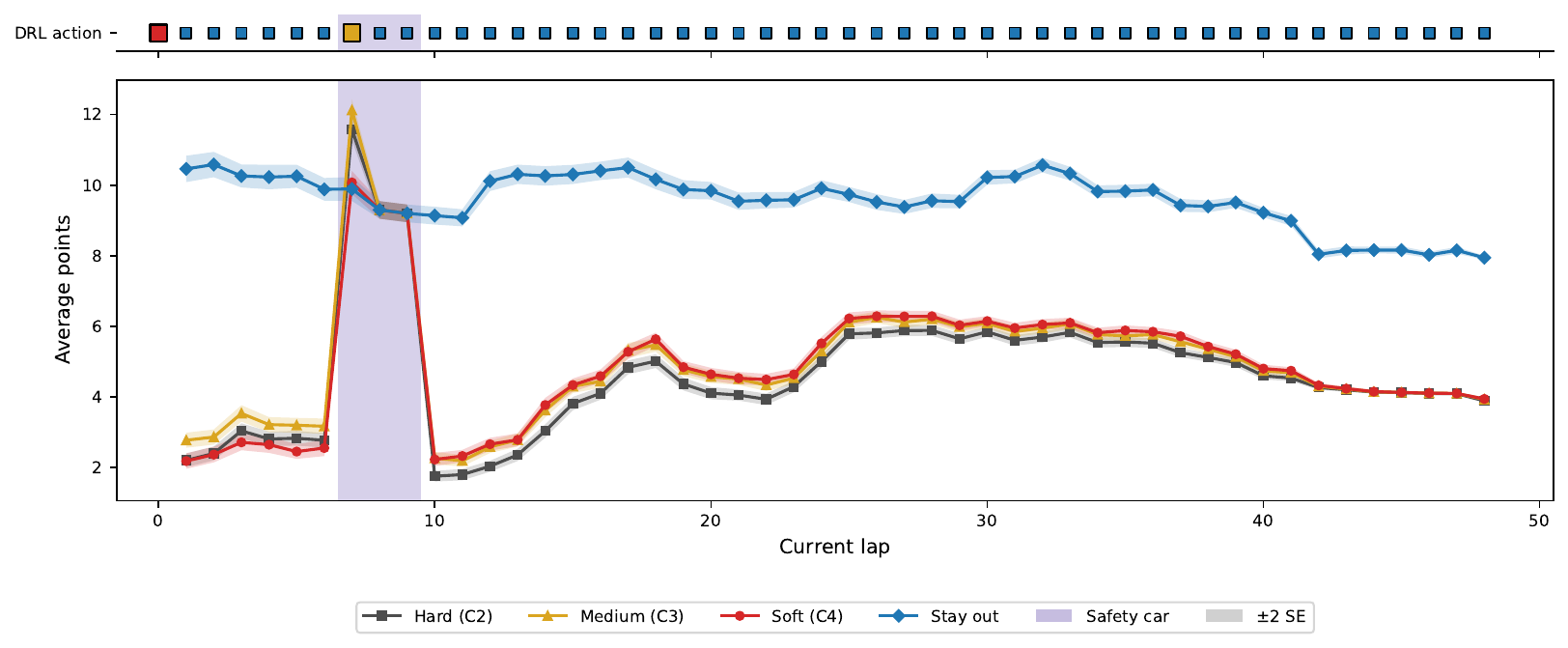}
\caption{Average championship points by lap for each candidate action over $M = 500$ rollouts. 2024 Saudi Arabian Grand Prix, PPO policy driving in place of Lando Norris.}
\label{fig:decision_evolution}
\end{figure}

Figure~\ref{fig:decision_evolution} shows the average championship points of each candidate action at every lap of the race, with a band of $\pm 2$ SE. The ``stay out'' action defers to the policy from lap $n+1$ onward, so it already includes the stop that the policy itself would make: at each lap the comparison is between stopping now and letting the policy choose when to stop. The top strip reports the action recommended by the policy at each lap.

\new{During the first six laps, pitting for any compound earns fewer expected points than staying on track, since the driver would abandon the initial soft-tire stint prematurely; consistent with this, the DRL policy recommends ``no pit'' throughout the opening phase. On lap $7$, when the safety car is deployed, pitting for medium tires becomes preferable to staying on track, with hard tires close behind, and this is the lap in which the policy prescribes its stop. On laps $8$ and $9$ the four actions coincide, since no stop is possible on them. From lap $10$ onward any further stop would be a second one, and staying on track remains the preferred continuation over the rest of the race.}

The look-ahead can also be applied retrospectively to review the decisions taken in 2024. For every lap $n$ of the eight races, the state $s_n$ replays the driver's actual starting compound and pit stops up to lap $n-1$, and the recommended action is compared with the actual one through the paired difference in points over the $M = 500$ rollouts. Since the actual starting compound often differs from the one the policy would choose, the analysis uses the PPO policy further trained with a randomly drawn starting compound (Appendix~\ref{app:realhist}). The recommendation coincides with the actual decision in $94\%$ of the laps, most of them laps in which neither stops. Of the $28$ laps where the two differ, the difference exceeds two standard errors in $12$, $7$ of them in favor of the recommendation (mean difference of $0.73$ points) and $5$ in favor of the actual decision ($0.23$ points) (Figure~\ref{fig:realhist} in Appendix~\ref{app:realhist}).

\begin{figure}[htbp]
\centering
\includegraphics[width=1.0\linewidth,trim=0 10bp 0 0,clip]{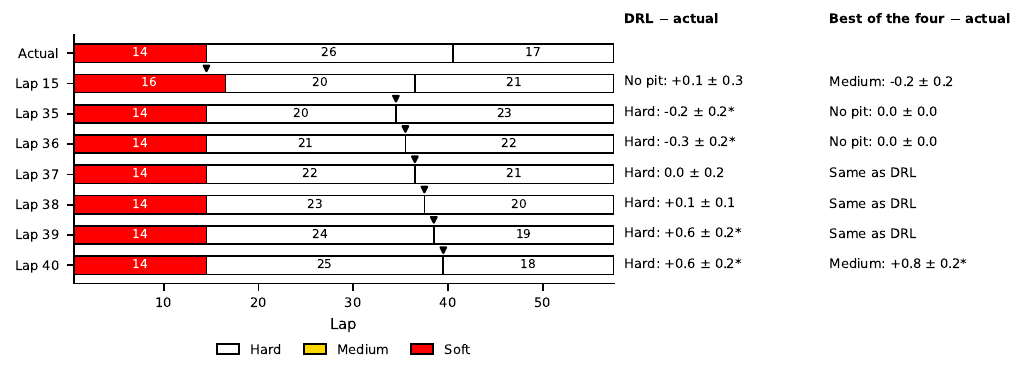}
\caption{Actual 2024 tire stints of Alonso in the Bahrain Grand Prix and stints implied by the DRL recommendation at every lap where it differs from the actual decision. Each lower bar replays the race up to the lap marked by the triangle and shows one representative continuation under the PPO policy, without further yellow flags. On the right, the mean paired difference in points against the actual decision for the recommended action and for the best of the four actions, $\pm 2$ SE over $M = 500$ rollouts with simulated yellow flags; positive values favor the recommended action (or the best of the four) over the actual decision, and $^{*}$ marks differences larger than 2 SE. The best action is selected with a separate set of $500$ rollouts, so its difference is not biased upward by the selection.}
\label{fig:realhist_best}
\end{figure}

The look-ahead evaluates all four actions, so it can be used to check the recommendation. Figure~\ref{fig:realhist_best} shows the laps of the Bahrain Grand Prix where the recommendation for Alonso differs from his actual decision, together with the best of the four actions at each of them. At lap $15$ Alonso stopped and the policy recommends staying out, but the rollouts show no clear difference between the two. Alonso's next stop was at lap $41$, whereas from lap $35$ on the policy recommends stopping for hard tires. At laps $35$ and $36$ the rollouts favor staying out, at laps $37$ and $38$ they show no clear difference, and at laps $39$ and $40$ they favor the recommendation, although at lap $40$ medium tires would be even better. The tire sets available at each decision point are reported in Table~\ref{tab:tyre_sets_decisions} of Appendix~\ref{app:tyre_sets}.

The rankings of the four actions inherit the biases of the policy used to roll out the remaining laps: where it is suboptimal, they may be distorted. PPO, one of the best-performing methods in our evaluation, mitigates this concern but does not eliminate it. Computation times, including those of the look-ahead, are reported in Appendix~\ref{app:runtimes}.

The trained policy is used as a simulator rather than as an autonomous agent: the strategist sees the outcomes of each action and decides how to combine them with information the model does not have.

\section{Conclusions}\label{sec:conclusions}

\rev{This paper studied Formula~1 pit-stop strategy as a sequential decision problem under uncertainty and competition, with the aim of supporting rather than replacing human strategists. We trained five DRL methods in a race simulation environment calibrated with historical data, which models lap times, drivers' recent performance, overtaking, rivals' pit-stop decisions, and yellow flag events, and benchmarked them against a QMIP formulation that is optimal in deterministic single-driver settings, where all DRL methods except TRPO achieve near-optimal performance.}

\new{With a full grid of competitors, PPO and Recurrent PPO attain the best average final positions in fully simulated races, and A2C and PPO when races are replayed with the rivals' decisions and yellow flags observed in 2024. PPO and Recurrent PPO are the only methods that significantly improve on the QMIP, on average over the eight circuits, in both settings, PPO by $0.33$ and $0.49$ positions in each setting and Recurrent PPO by less, while the QMIP, which cannot anticipate yellow flag events, remains competitive with the other DRL methods. TRPO is among the weakest methods in every setting considered, so the choice of algorithm still matters.}

\rev{We then used the trained PPO policy as a look-ahead simulator that reports the distribution of final positions and championship points under each candidate action. \new{In an illustrative case from the 2024 Saudi Arabian Grand Prix, on the lap in which a safety car is deployed, it shows that stopping for medium or hard tires is preferable to both staying out and a second soft-tire stint, and it identifies the laps at which deviating from the policy is costly.} Applied retrospectively to eight 2024 races, the look-ahead identifies the laps where the actual decision, the recommendation, or a third action was preferable, so it can be used to check the recommendation itself. A strategist can weigh the recommendation against information the model lacks.}

\rev{Model parameters are estimated from pre-race information only; with live data, the QMIP could be re-solved with updated parameters during the race, whereas the DRL policies would need retraining. Rivals follow data-driven behavioral models rather than responding strategically to the agent, a limitation that multi-agent formulations could address, also allowing for coordination between teammates.
\new{The tire sets available at the start of the race are inferred from the sets used in each session and the 2024 Sporting Regulations, since the returned sets are not reported, so the availability of tires at each decision point is approximate.}
Finally, evaluating the look-ahead tool with practicing strategists would show how it affects their reliance on the policy and the quality of their decisions.}

\backmatter

\bmhead{Acknowledgements}

This work was supported by the Agencia Nacional de Investigaci\'on y Desarrollo (ANID), Fondecyt Iniciaci\'on 11241531. Powered@NLHPC: This research was partially supported by the supercomputing infrastructure of the NLHPC (ECM-02).

\bmhead{Use of AI tools}

The authors used Claude (Anthropic) to assist with code development, data analysis and drafting of parts of the text. The authors reviewed all outputs and take full responsibility for the content of this article.

\bmhead{Statements and Declarations}

\noindent\textbf{Funding.} This work was funded by the Agencia Nacional de Investigaci\'on y Desarrollo (ANID), Fondecyt Iniciaci\'on 11241531.

\noindent\textbf{Competing interests.} The authors declare no competing interests.

\noindent\textbf{Ethics approval and consent to participate.} Not applicable.

\noindent\textbf{Consent for publication.} Not applicable.

\noindent\textbf{Data availability.} Race data were obtained from the public Formula~1 timing feeds through FastF1 \citep{fastf1} and \cite{openf1}. The processed data are available at \cite{mendez_thraves_2025_data}.

\noindent\textbf{Materials availability.} Not applicable.

\noindent\textbf{Code availability.} The code and trained models are available at \cite{mendez_thraves_2025_data}.

\noindent\textbf{Author contribution.} F.M. contributed to data curation, formal analysis, investigation, methodology, software, visualization, and writing (review and editing). C.T. contributed to conceptualization, formal analysis, funding acquisition, methodology, software, supervision, and writing (original draft).

\bigskip

\bibliography{sn-bibliography}

\clearpage
\begin{appendices}
\numberwithin{figure}{section}\numberwithin{table}{section}
\numberwithin{equation}{section}\numberwithin{algorithm}{section}
\renewcommand{\thefigure}{\thesection.\arabic{figure}}
\renewcommand{\thetable}{\thesection.\arabic{table}}
\renewcommand{\theequation}{\thesection.\arabic{equation}}
\renewcommand{\thealgorithm}{\thesection.\arabic{algorithm}}

\FloatBarrier
\new{\section{Rescaling of tire age in 2021}\label{app:tirelife2021}
The tires used from 2022 onward have a different construction from those of 2021, so the lap-time loss per lap of tire age differs between the two periods. We express the 2021 tire age in equivalent laps of 2022--2023 by multiplying it by a factor $f$, common to all compounds and circuits, and estimate $f$ in two ways.

The first estimate uses lap times. We estimate Equation~\eqref{eq:reg} for each circuit with 2021--2023 data, imposing that the 2021 tire age slope of each compound equals $f$ times its 2022--2023 slope, with $f$ common to all circuits. Since fuel consumption and track evolution shape the lap-time trend over the race differently across seasons, the lap number terms are specific to each year. Using the 15 circuits with laps in both periods (15,524 laps in 2021), we obtain $f=0.69$, with a standard error of 0.069 clustered by race.

The second estimate uses pit-stop decisions. We estimate a logit model of the decision to pit on 2021--2023 laps, with indicators for SC and VSC and for having used two compounds, compound-specific tire age slopes, the driver's position, circuit effects, a 2021 effect, and a function of the fraction of the race completed, specific to 2021 (a B-spline with five degrees of freedom). As in the rivals' model (Section~\ref{sec:environment_rivals}), each SC or VSC period contributes one observation per driver, on its first lap, equal to one if the driver stopped at any lap of the period. With the 2021 tire age multiplied by $f$, the likelihood is maximized at $f=0.87$, with a standard error of 0.074 clustered by race.

The two estimates are compatible ($p=0.08$). Treating them as independent, since they rely on different data, we combine them by a precision-weighted average, which gives $f=0.78$ (95\% confidence interval 0.68--0.88).}

\FloatBarrier
\new{\section{Joint model for tire coefficients}\label{app:regTC_tirecoef_matrix_completion}

Estimating Equation~\eqref{eq:reg} separately for each circuit yields compound coefficients that are often imprecise and frequently violate the expected ordering of compounds. This occurs when a compound was used in few stints at a circuit, when it was not used at all in 2021--2023, or when tire age is highly correlated with the lap number, as in the Saudi Arabian GP, where pit stops were synchronized by a red flag and SC periods in 2021, and by an SC period in 2022 and 2023. We therefore replace the compound coefficients of Equation~\eqref{eq:reg} by those of a model fitted jointly across circuits.

Let $o_{gc}=\beta_{c,1}$ and $s_{gc}=\beta_{c,2}$ denote the intercept and degradation coefficients of compound $c$ at circuit $g$. The joint model is
\begin{equation*}
o_{gc}=\lambda_g\, a_c, \qquad s_{gc}=\exp\left(\theta_g+\lambda'_g\, b_c\right),
\end{equation*}
where $a_c$ and $b_c$ are compound profiles common to all circuits, with $a_c$ decreasing and $b_c$ increasing in the softness of the compound, $\lambda_g,\lambda'_g>0$ are circuit-specific scale factors, and $\theta_g$ is the degradation level of circuit $g$. The ordering is imposed by writing each profile as a cumulative sum of positive steps between consecutive compounds; since the 2021--2023 data are not informative about the step between C4 and C5, it is set equal to the step between C3 and C4. The level $\theta_g$ is a random effect, $\theta_g\sim N(\mu,\tau^2)$, with $\tau=0.768$ estimated by maximum likelihood, which shrinks the degradation level of circuits with poorly identified slopes toward the common mean. The scale factors are penalized toward one ($\log\lambda_g$ and $\log\lambda'_g$ with standard deviation $0.3$). The model is fitted by generalized least squares to the per-circuit estimates of Equation~\eqref{eq:reg} for the 21 circuits of the 2024 calendar with 2021--2023 data, each weighted by the inverse of its estimated covariance matrix. Since the intercepts of Equation~\eqref{eq:reg} are relative to a reference compound, only their differences within each circuit enter the fit. Table~\ref{tab:tire_coefs} reports the resulting coefficients for the eight circuits considered.

\begin{table}[htbp]
\centering
\begin{tabular}{llccc}
\toprule
Grand Prix & Compound & Stints 2021--2023 & $\beta_{c,1}$ & $\beta_{c,2}$ \\
\midrule
Bahrain & C1 & 33 & 0.000 & 0.1171 \\
 & C2 & 45 & -0.206 & 0.1388 \\
 & C3 & 109 & -0.415 & 0.1674 \\
\midrule
Belgian & C2 & 16 & 0.000 & 0.1574 \\
 & C3 & 60 & -0.145 & 0.1780 \\
 & C4 & 33 & -0.256 & 0.1919 \\
\midrule
Dutch & C1 & 0 & 0.000 & 0.0408 \\
 & C2 & 41 & -0.188 & 0.0486 \\
 & C3 & 50 & -0.378 & 0.0589 \\
\midrule
Hungarian & C3 & 75 & 0.000 & 0.0574 \\
 & C4 & 45 & -0.178 & 0.0661 \\
 & C5 & 3 & -0.356 & 0.0761 \\
\midrule
Miami & C2 & 39 & 0.000 & 0.0163 \\
 & C3 & 40 & -0.123 & 0.0198 \\
 & C4 & 5 & -0.217 & 0.0224 \\
\midrule
Saudi Arabian & C2 & 50 & 0.000 & 0.0111 \\
 & C3 & 62 & -0.167 & 0.0161 \\
 & C4 & 5 & -0.295 & 0.0202 \\
\midrule
Singapore & C3 & 18 & 0.000 & 0.0405 \\
 & C4 & 21 & -0.204 & 0.0465 \\
 & C5 & 3 & -0.409 & 0.0535 \\
\midrule
United States & C2 & 79 & 0.000 & 0.0782 \\
 & C3 & 86 & -0.156 & 0.0885 \\
 & C4 & 5 & -0.275 & 0.0955 \\
\bottomrule
\end{tabular}
\caption{Tire coefficients of the joint model for each Grand Prix and compound. $\beta_{c,1}$ is expressed relative to the hardest compound of the Grand Prix (seconds) and $\beta_{c,2}$ is the lap-time increase per lap of tire age (seconds per lap). The number of stints run with each compound in 2021--2023 indicates how much each coefficient relies on the circuit's own data.}
\label{tab:tire_coefs}
\end{table}}

\FloatBarrier
\new{\section{Specifications for updating driver coefficients}\label{app:adjust_driver}

We compare four specifications of Equation~\eqref{eq:regcoef} that differ in how lap times are expressed: (i) in seconds, with the pole lap time as a proxy for circuit length; (ii) as differences to the reference driver; (iii) as differences to the session median; and (iv) relative to lap time, as in Equation~\eqref{eq:regcoef}. Each specification is estimated using only the races prior to the target race, and evaluated on every race of 2023 and 2024 against the coefficient estimated with Equation~\eqref{eq:reg} on that race alone. Two benchmarks without correction are included: the driver's coefficient in the previous race, and the coefficient estimated with Equation~\eqref{eq:reg} on previous editions of the same circuit. Table~\ref{tab:reg_coef_drivers} reports the mean squared errors. All specifications substantially reduce the error relative to both benchmarks, and the relative specification attains the lowest error in both seasons. Dropping $X_{\text{Q}}$ (the driver's best qualifying lap time relative to that of the reference driver) changes the error by at most $0.02$, whereas using fewer previous races increases it, up to about $0.35$ when none is available. Accordingly, $X_{\text{Q}}$ is dropped when the qualifying session is missing or was run in wet conditions, and drivers with fewer than three previous races use the available ones. Races in which a driver, or the reference driver, completes fewer than 15 clean laps are not used, either as observations or as previous races.

\begin{table}[htbp]
\centering
\begin{tabular}{lcccccc}
\toprule
 & \multicolumn{4}{c}{Specification} & \multicolumn{2}{c}{No correction} \\ \cmidrule(lr){2-5}\cmidrule(lr){6-7}
Season & (i) & (ii) & (iii) & (iv) & Previous race & Same circuit \\ \midrule
2023 & 0.136 & 0.153 & 0.166 & 0.135 & 0.249 & 0.538 \\
2024 & 0.221 & 0.195 & 0.210 & 0.191 & 0.322 & 0.618 \\ \bottomrule
\end{tabular}
\caption{Out-of-sample mean squared error, in (seconds per lap)$^2$, of the predicted driver coefficients for the races of 2023 (189 driver--race pairs) and 2024 (276), by specification of Equation~\eqref{eq:regcoef} and for the two benchmarks without correction.}
\label{tab:reg_coef_drivers}
\end{table}}

\FloatBarrier
\new{\section{Lap times under yellow flag conditions}\label{app:regLTYF}

Equation~\eqref{eq:reg_yf} is estimated with a linear mixed model on lap times divided by $L^{(g)}$, using the 2021--2023 laps that a driver ran entirely under VSC or SC (1117 laps from 34 races). Laps in which a car went through the pit lane without changing tires, as when the field is led through the pit lane behind the SC, are excluded, also from the excess times below. The reference lap time $L^{(g)}$ is the median clear-track lap without pit stops at circuit $g$, and $G^{(g)}_{\text{in}}$ and $G^{(g)}_{\text{out}}$ are the median excess times of in-laps and out-laps over the cars that did not pit in the same lap. We exclude the lap in which a SC is deployed, since each car runs it only partially under SC, and pit stops made once the field is queued, whose duration mixes the regrouping modeled in Algorithm~\ref{alg:simevolution}. The model includes a random effect common to all cars in each race lap and the circuit effect $u_g$ under SC.

Table~\ref{tab:yf_common} reports the common coefficients. A lap under VSC is about $34\%$ slower than a clear-track lap, and a lap under SC about $42\%$ slower before the circuit effect. The pace of the first SC lap, before the field is queued, does not differ significantly ($1.432$, s.e. $0.027$), so a single SC coefficient is used. The pit-lane factors $\kappa_{\text{in}}$ and $\kappa_{\text{out}}$ are common to SC and VSC. Table~\ref{tab:yf_gp} reports the circuit-level quantities and the resulting lap times and pit-lane losses. Across the eight circuits, a pit stop under yellow flag conditions costs between $59\%$ and $69\%$ of the clear-track loss, which saves between $6.3$ and $12.3$ seconds.

\begin{table}[htbp]
\centering
\begin{tabular}{lcc}
\toprule
Coefficient & Estimate & Std. Err. \\
\midrule
$\rho_{\text{VSC}}$ & 1.337 & (0.023) \\
$\rho_{\text{SC}}$ & 1.420 & (0.026) \\
$\kappa_{\text{in}}$ & 0.366 & (0.227) \\
$\kappa_{\text{out}}$ & 0.702 & (0.092) \\
Std. dev. of $u_g$ & 0.087 & \\
\bottomrule
\end{tabular}
\caption{Common coefficients of the yellow flag lap time model (Equation~\eqref{eq:reg_yf}), estimated with 2021--2023 data.}
\label{tab:yf_common}
\end{table}

\begin{table}[htbp]
\centering
\small
\setlength{\tabcolsep}{4pt}
\begin{tabular}{lcccccccc}
\toprule
Grand Prix & $L^{(g)}$ & $G^{(g)}_{\text{in}}$ & $G^{(g)}_{\text{out}}$ & $u_g$ & VSC lap & SC lap & Pit (YF) & Pit (clear) \\
\midrule
Bahrain & 98.315 & 3.227 & 20.841 & -0.129 & 131.5 & 127.0 & 15.8 & 24.1 \\
Saudi Arabian & 95.170 & 3.347 & 19.085 & 0.033 & 127.3 & 138.3 & 14.6 & 22.4 \\
Miami & 92.768 & 5.034 & 14.860 & -0.076 & 124.1 & 124.7 & 12.3 & 19.9 \\
Hungarian & 84.573 & 3.069 & 18.674 & 0.000 & 113.1 & 120.1 & 14.2 & 21.7 \\
Belgian & 112.776 & 5.122 & 13.653 & 0.000 & 150.8 & 160.2 & 11.5 & 18.8 \\
Dutch & 76.390 & 3.681 & 18.703 & -0.092 & 102.2 & 101.5 & 14.5 & 22.4 \\
Singapore & 100.414 & 9.947 & 19.977 & -0.023 & 134.3 & 140.4 & 17.7 & 29.9 \\
United States & 102.338 & 0.678 & 19.619 & -0.026 & 136.9 & 142.7 & 14.0 & 20.3 \\
\bottomrule
\end{tabular}
\caption{Circuit-level quantities of the yellow flag lap time model and resulting lap times and pit-lane losses, in seconds. Pit (YF) is the pit-lane loss $\kappa_{\text{in}}G^{(g)}_{\text{in}}+\kappa_{\text{out}}G^{(g)}_{\text{out}}$ and Pit (clear) is $G^{(g)}_{\text{in}}+G^{(g)}_{\text{out}}$. No SC laps remain in the estimation sample at the Hungarian and Belgian GPs in 2021--2023, so their circuit effect is set to zero.}
\label{tab:yf_gp}
\end{table}}

\FloatBarrier
\new{\section{Tire Sets Available at the Start of the Race}\label{app:tyre_sets}
The tire sets used in each session of the 2024 weekends (practice, sprint qualifying, sprint and qualifying) and their ages are taken from FastF1 and OpenF1 (\url{https://openf1.org}). Under Article 30.2 of the 2024 Sporting Regulations, a driver receives 13 dry sets in a regular weekend and returns two after each practice session and, if they reach Q3, one soft set after qualifying; in sprint weekends, 12 sets are allocated and one is returned after FP1, one after the sprint and three after qualifying. One hard and one medium set, the race specifications at the eight circuits, cannot be returned. Which set is returned is not reported; we assume that, among the sets not used later in the weekend, teams return the one with the most laps, preferring soft sets. The resulting inventory has 6.73 sets per driver on average, 55\% of them new (counting only push laps). In 2024, 17\% of the race stints started on a used set (15\% of the starts and 19\% of the stints after a pit stop), and 192 of these 196 stints match a set of the inferred inventory. Table~\ref{tab:tyre_sets_decisions} reports the sets available to each driver at the decision points reviewed in Figure~\ref{fig:realhist}, which bound the length of the stints shown there.}

\begin{table}[htbp]
\centering
\resizebox{\textwidth}{!}{%
\begin{tabular}{llcclll}
\toprule
 & & Laps & On car & \multicolumn{3}{c}{Sets available (age in laps)} \\
\cmidrule(lr){5-7}
GP (driver) & Lap & left & (age) & Hard & Medium & Soft \\
\midrule
Bahrain (Alonso) & 15 & 43 & C3 (15) & C1: 0, 0 & C2: 0 & C3: 1, 1 \\
 & 35--40 (39) & 19 & C1 (23) & C1: 0 & C2: 0 & C3: 1, 1, 16 \\
\addlinespace
Saudi Arabian (Norris) & 7 & 44 & C3 (6) & C2: 0 & C3: 0 & C4: 1, 3, 3 \\
 & 29 & 22 & C3 (25) & C2: 0 & C3: 0 & C4: 1, 3, 3 \\
 & 31--36 (34) & 17 & C3 (30) & C2: 0 & C3: 0 & C4: 1, 3, 3 \\
\addlinespace
Miami (Ocon) & 22 & 36 & C3 (21) & C2: 0 & C3: 2, 2 & C4: 1, 1, 1 \\
\addlinespace
Hungarian (Tsunoda) & 29 & 42 & C4 (28) & C3: 0, 0 & C4: -- & C5: 1, 1, 2 \\
\addlinespace
Belgian (Gasly) & 9 & 36 & C3 (8) & C2: 0, 0 & C3: -- & C4: 0, 0, 0, 0 \\
 & 25 & 20 & C2 (15) & C2: 0 & C3: 9 & C4: 0, 0, 0, 0 \\
 & 26--28 (26) & 19 & C2 (16) & C2: 0 & C3: 9 & C4: 0, 0, 0, 0 \\
\addlinespace
Dutch (Stroll) & 30 & 43 & C2 (29) & C1: 0 & C2: -- & C3: 1, 1, 1, 2 \\
\addlinespace
Singapore (Hülkenberg) & 29 & 34 & C4 (28) & C3: 0 & C4: -- & C5: 1, 1, 1, 1 \\
\addlinespace
United States (Magnussen) & 3 & 54 & C3 (2) & C2: 0 & C3: 0, 2 & C4: 1, 1, 1 \\
 & 17 & 40 & C3 (13) & C2: 0 & C3: 0, 2 & C4: 1, 1, 1 \\
 & 36--37 (36) & 21 & C2 (18) & C2: -- & C3: 0, 2, 14 & C4: 1, 1, 1 \\
\bottomrule
\end{tabular}}
\caption{Tire sets available to the driver at each decision point of Figure~\ref{fig:realhist}, with the age in laps of every set of each compound and of the set on the car. For a group of consecutive laps, the lap of the largest difference is given in parentheses and the state refers to that lap. A pit stop fits the least-worn set of the chosen compound. A dash means that no spare set of that compound is left; the set taken off returns to the driver's sets with its age, so the set on the car can always be refitted, with the age it has at that lap. Ages follow Section~\ref{sec:environment}: laps under SC, VSC or red flag, and the non-push laps of earlier sessions, do not count.}\label{tab:tyre_sets_decisions}
\end{table}

\FloatBarrier
\section{Safety Car and Virtual Safety Car Event Data}\label{app:prob_yf}
Table~\ref{tab:sc_vsc_probs} reports, for each circuit, the estimated probability that at least one SC or at least one VSC occurs during a race. In the simulation, the probability of an event in each lap is $1-(1-P)^{1/N}$, where $P$ is the probability in the table and $N$ is the number of laps of the race. The SC probability of the Saudi Arabian GP is capped at $0.70$, since a probability of one would imply a SC in every lap of the simulation. The last two columns give the laps with SC or VSC conditions in the actual 2024 race, where a dash means that none occurred; these laps are used when the yellow flag events are taken from the real race.

\begin{table}[htbp]
\centering
\begin{tabular}{lcccc}
\toprule
& \multicolumn{2}{c}{Estimated Probability} & \multicolumn{2}{c}{Actual Laps (2024 GP)} \\
\cmidrule(lr){2-3}\cmidrule(lr){4-5}
Grand Prix & SC & VSC & SC & VSC \\
\midrule
Bahrain         & 0.37 & 0.22 & ---    & ---    \\
Saudi Arabian   & 0.70 & 0.67 & \new{7--9}  & ---    \\
Miami           & \new{0.50} & \new{0.50} & \new{28--32} & 22--23 \\
Hungarian       & 0.25 & 0.33 & ---    & \new{---}   \\
Belgian         & 0.64 & 0.11 & ---    & ---    \\
Dutch           & 0.67 & 0.67 & ---    & ---    \\
Singapore       & 0.93 & 0.43 & ---    & ---    \\
United States   & 0.36 & 0.38 & \new{3--5}    & ---    \\
\bottomrule
\end{tabular}
\caption{Estimated probability of at least one SC or VSC event per race, \new{taken from Aguad and Thraves (2024) except for Miami, which is not included there and is estimated from its 2022--2023 races,} and lap ranges during which SC or VSC conditions were active in the corresponding 2024 race.}
\label{tab:sc_vsc_probs}
\end{table}

\FloatBarrier
\section{Race Evolution}\label{app:race_evolution}

Algorithm~\ref{alg:simevolution} describes the evolution of one lap of the simulation, given the state at its beginning and the agent's action.

\begin{algorithm}[htbp]
\caption{System evolution for a lap}\label{alg:simevolution}
\begin{algorithmic}[1]
\Require{Lap $n$, state at the beginning of the lap $S$, agent's action}
\State \rev{$\ell \gets$ number of consecutive laps (including $n$) with $S_{SC}\geq 1$; $\ell=0$ if $S_{SC}=0$}
\State Sample rivals' pit-stop decisions using LogitStop (pit stops are only allowed on clear-track laps and on the first lap of an SC or VSC period; otherwise the agent's pit action is a no-op) 
\State For rival drivers who stop, get the new tire compound with CondLogitCompound 
\If{$S_{SC}+S_{VSC}=0$} \Comment{Clear track lap}
    \State Compute $\{\hat{y}_{d}\}_{d\in\mathcal{D}}$ with RegLTC \Comment{Preliminary lap times; each driver's coefficient includes the prediction error drawn at the start of the race}
\Else \Comment{There is a SC or VSC}
    \State Compute $\{\hat{y}_{d}\}_{d\in\mathcal{D}}$ with RegLTYF \Comment{Preliminary lap times\rev{; the SC/VSC coefficient is used from the first lap of the period}}
\EndIf
\State \rev{\textbf{if} $n$ is the final lap \textbf{and} driver $d$ pits \textbf{then} charge both PitIn and PitOut to $\hat{y}_{d}$} \Comment{\rev{No out-lap exists after the final lap}}
\If{$S_{SC}\geq 1$ \textbf{and} $\ell \geq 3$}
    \State Queue the field behind the leader, with gaps of $t_{\text{SC\_gap}}$ between consecutive cars; positions are frozen; \textbf{return} \Comment{No overtaking}
\EndIf
\State $d^{(i)} \gets$ driver in $i$-th position at the beginning of lap $n$
\For{$h = 2$ to $20$}
    \State $d^{(h)} \gets$ driver in the $h$\textsuperscript{th} position at the beginning of lap $n$
    \For{$k = h$ to $2$}
        \State $d^{(k-1)} \gets$ driver (partially) in the $k-1$\textsuperscript{th} position 
        \State $X_{\Delta_{2}}^{(k,k-1)} \gets$ time gap at the beginning of lap $n$ between drivers $d^{(k)}$ and $d^{(k-1)}$
        \State $X_{\text{PredGap}} \gets X_{\Delta_{2}}^{(k,k-1)} + \hat{y}_{d^{(k)}} - \hat{y}_{d^{(k-1)}}$
        \If{$S_{SC}+S_{VSC}=0$} \Comment{Clear track lap}
            \State Get $\tau\leftarrow\text{LogitOvertake}$  \Comment{Probability that $d^{(k)}$ overtakes $d^{(k-1)}$, evaluated with the mean lap times, since the logit already contains their uncertainty}
            \State Overtake $\gets \text{Uniform}(0,1) < \tau$
        \Else \Comment{There is a SC or VSC}
            \State Overtake $\gets X_{\text{PredGap}} < 0$
        \EndIf
        \If{Overtake} \Comment{Overtake is predicted to occur}
            \State $\hat{y}_{d^{(k)}} \gets \max\{\hat{y}_{d^{(k)}}, \hat{y}_{d^{(k-2)}}-  X_{\Delta_{2}}^{(k,k-2)} + t_{\text{min\_gap}}\}$
            \State $\hat{y}_{d^{(k-1)}} \gets \max\{\hat{y}_{d^{(k-1)}}, \hat{y}_{d^{(k)}} + X_{\Delta_{2}}^{(k,k-1)}+ t_{\text{min\_gap}}\}$
            \State Swap $d^{(k)}$ and $d^{(k-1)}$
        \Else \Comment{Overtake does not occur}
            \State $\varepsilon \sim \text{Exp}(0.26)$, with mean $0.26$~s \Comment{\rev{Following gap: $t_{\text{min\_gap}}+\varepsilon$ seconds, resampled per pair and lap}}
            \State $\hat{y}_{d^{(k)}} \gets \max\{\hat{y}_{d^{(k)}}, \hat{y}_{d^{(k-1)}} - X_{\Delta_{2}}^{(k,k-1)} + t_{\text{min\_gap}} \rev{+ \varepsilon}\}$
            \State \textbf{break}
        \EndIf
    \EndFor
\EndFor
\If{$S_{SC}\geq 1$ \textbf{and} $\ell = 2$}
    \State Scale the gaps to the leader by $\tfrac{1}{2}$, after applying the lap times \Comment{Progressive regrouping behind the SC}
\EndIf
\end{algorithmic}
\end{algorithm}

\rev{Three remarks complement Algorithm~\ref{alg:simevolution}. First, during a SC period the field is regrouped progressively: gaps are halved on the second lap of the period and, from the third lap on, drivers queue behind the leader at $t_{\text{SC\_gap}}=1.5$~s (the median gap observed in real SC queues), with track positions frozen; pit stops are only possible on the first lap of the period, when teams react to the deployment. \new{This is a simplification, since the pit lane remains open during the whole period; in practice, however, most stops under a safety car are made on its first lap (about two thirds in the dry races of 2021--2024), and in a sensitivity check that also allows the agent to stop on the second lap, such a stop yields a similar expected outcome to one on the first lap, so the recommendation to stop under the safety car in Figure~\ref{fig:decision_support} is unchanged.} The same restriction applies under VSC, to the agent and to the rivals: a stop is only possible on the first lap of a VSC period. Second, the following gap of a car that fails to overtake is not a constant: it is drawn as $t_{\text{min\_gap}}+\varepsilon$ with $t_{\text{min\_gap}}=0.4$~s and $\varepsilon$ exponential with mean $0.26$~s, so that the mean following gap matches the $0.66$~s observed for cars that are faster than the car ahead but fail to pass it. A constant gap would keep groups of cars at a fixed spacing. Third, a pit stop on the final lap is charged the full pit-lane time (PitIn and PitOut) within that lap, since no out-lap exists; otherwise such a stop would cost only the pit-entry component.}

\FloatBarrier
\section{Models Architecture and Seed Selection}\label{app:architecture}

\subsection{Hyperparameters}

Tables~\ref{tab:params_no_competition} and~\ref{tab:params_all_rivals} report the hyperparameters of the five DRL algorithms in the simplified setting without competition and in the setting with a full grid of competitors, respectively. The algorithms are those of Stable-Baselines3 and its contrib package (TRPO and Recurrent PPO), and the parameters that are not listed take their default values. Without competition, each configuration is trained with five seeds.

\begin{table}[htbp]
\centering
\setlength{\tabcolsep}{4pt}
\renewcommand{\arraystretch}{0.92}
\new{\begin{tabular}{lccccc}
\toprule
Hyperparameter & Double DQN & A2C & TRPO & PPO & Recurrent PPO \\
\midrule
Learning rate               & $1 \!\times\! 10^{-4}$ & $3 \!\times\! 10^{-4}$ & $1 \!\times\! 10^{-3}$ & $3 \!\times\! 10^{-4}$ & $5 \!\times\! 10^{-4}$ \\
Discount factor ($\gamma$)  & \multicolumn{5}{c}{0.99} \\
\midrule
LSTM layers $\times$ units  & --            & --          & --            & --            & $1 \times 32$ \\
Hidden layers               & [64, 64, 64]  & [64, 64]    & [64, 64, 64]  & [64, 64, 64]  & [32, 32]      \\
\midrule
$n_\mathrm{steps}$          & --            & $n_l$       & $n_l$         & $2 n_l$       & $n_l$         \\
Batch size                  & $2 n_l n_e$   & --          & $n_l n_e$     & $n_l n_e$     & $n_l n_e$     \\
Entropy coefficient         & --            & 0.05        & --            & 0.01          & 0.10          \\
GAE $\lambda$               & --            & 0.95        & 0.95          & 0.95          & 0.95          \\
Clip range                  & --            & --          & --            & 0.20          & 0.30          \\
Update epochs               & --            & --          & --            & 10            & 10            \\
\bottomrule
\end{tabular}}
\caption{Hyperparameters of each RL algorithm without competition. $n_l$: laps per circuit; $n_e = 10$: parallel environments. All use ReLU activation, reward normalization and 5M timesteps. Double DQN uses a replay buffer of $3 \times 10^6$; A2C uses RMSProp with normalized advantages; PPO uses target KL $= 0.015$; the LSTM of Recurrent PPO is used by the actor only.}
\label{tab:params_no_competition}
\end{table}

  \begin{table}[htbp]
  \centering
  \setlength{\tabcolsep}{4pt}
  \renewcommand{\arraystretch}{0.92}
  \begin{tabular}{lccccc}
  \toprule
  Hyperparameter & Double DQN & A2C & TRPO & PPO & Recurrent PPO \\
  \midrule
  Learning rate$^{\dagger}$   & $3 \!\times\! 10^{-4}$ & $5 \!\times\! 10^{-4}$ & $5 \!\times\! 10^{-4}$ & $3 \!\times\! 10^{-4}$ & $5 \!\times\!
  10^{-4}$ \\
  Discount factor ($\gamma$)  & \multicolumn{5}{c}{0.995} \\
  Reward normalization        & No & No & No & Yes & Yes \\
  \midrule
  LSTM layers $\times$ units  & --                & --          & --          & --          & $1 \times 64$ \\
  Hidden layers               & [256, 256, 256]   & [256, 128]  & [256, 128]  & [256, 128]  & [128]         \\
  \midrule
  $n_\mathrm{steps}$          & --          & $2 n_l$     & $2 n_l$     & $2 n_l$     & $2 n_l$       \\
  Batch size                  & 256         & --          & $2 n_l n_e$ & $2 n_l n_e$ & $2 n_l n_e$   \\
  Entropy coefficient$^{\ddagger}$ & --     & 0.01        & --          & 0.05  & 0.05          \\
  VF coefficient              & --          & 0.50        & --          & 0.01        & 0.01          \\
  GAE $\lambda$               & --          & 1.00        & 0.95        & 0.98        & 0.95          \\
  Clip range                  & --          & --          & --          & 0.20        & 0.25          \\
  Update epochs               & --          & --          & --          & 3           & 4             \\
  Target KL                   & --          & --          & 0.02        & --          & --            \\
  \bottomrule
  \end{tabular}
  \caption{Hyperparameters of each RL algorithm with competition. $n_l$: laps per circuit; $n_e = 10$: parallel environments. All use ReLU activation and 5M timesteps with early stopping (patience of 25 evaluations, 50 episodes each). Double DQN uses a replay buffer of $1.5 \!\times\! 10^{6}$, learning starts after 30{,}000 steps, exploration fraction 0.30 with $\epsilon_\mathrm{final} = 0.02$, and target network updates every 1{,}500 steps. $^{\dagger}$~Double DQN anneals the learning rate linearly from $3 \!\times\! 10^{-4}$ to $1 \!\times\! 10^{-5}$; the other algorithms use a constant rate. $^{\ddagger}$~With 0.02, PPO and Recurrent PPO collapsed their entropy early in training and kept stopping on the sixth lap of the Miami Grand Prix.}
  \label{tab:params_all_rivals}
  \end{table}

\FloatBarrier
\subsection{Seed Selection}

\begin{figure}[htbp]
\centering
\includegraphics[width=1.0\linewidth]{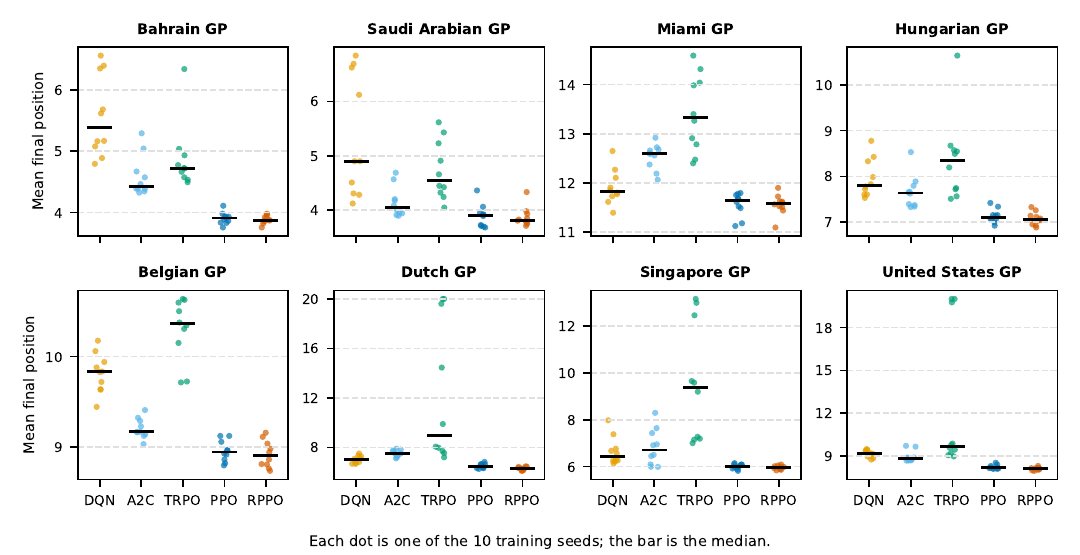}
\caption{Mean final position of the ten training seeds of each algorithm, by circuit, in the validation races used for seed selection.}
\label{fig:seed_spread}
\end{figure}

With competition, each agent is trained with ten seeds, and each seed is evaluated on $200$ simulated races independent of those of Figure~\ref{fig:stats_drivers_stochastic_2024}. A race is valid if the agent makes one to four pit stops, the first after lap 5 unless a safety car or virtual safety car has already been deployed, and, except on safety car or virtual safety car laps, none with seven or fewer laps to go and none fewer than eight laps after the previous stop. Among the seeds with at least $95\%$ valid races and no race disqualified for not using two compounds, the one with the best average final position is kept; if no seed meets both conditions, the one with the most valid races is kept. Only TRPO needs this fallback, in four of the eight circuits, where the selected seed reaches the share of valid races but leaves up to $3.5\%$ of them disqualified. Figure~\ref{fig:seed_spread} shows the mean final position of every seed in these validation races. The ten seeds of PPO and Recurrent PPO lie close together in every circuit. Double DQN and A2C vary more across seeds in some circuits, and TRPO is the most sensitive, with seeds far from the rest in several circuits. The agent evaluated in Section~\ref{sec:agent_competition} is the selected seed, which by construction is among the best of its method.

\FloatBarrier
\new{\section{Estimated Coefficients for the Clear-Track Lap Time Model}\label{app:table_reg_clean}

Table~\ref{tab:regression_cleanlap} reports the estimated coefficients of Equation~\eqref{eq:reg} for the 2024 United States Grand Prix, estimated with the 2021--2023 clear-track laps, before and after removing outliers, that is, laps at least 15 seconds slower than the median lap of their race (except in-laps, out-laps and the first lap) and the observations excluded by Algorithm~\ref{alg:IQR}. In the table, Interval\_behind and Interval\_front are $e^{-X_{\Delta_1}}$ and $e^{-X_{\Delta_2}}$, FirstLap\_pos is $X_{\text{1pos}}$, and TireLife is the effective tire age, with the 2021 laps rescaled as in Appendix~\ref{app:tirelife2021}. Removing $1.7\%$ of the laps increases the coefficient of determination from 0.877 to 0.966 and roughly halves the standard errors. Some coefficients also take the expected sign. For example, starting further back increases the first lap time, which holds only after the removal of outliers. Similarly, the tire life coefficient of the softest compound (C4) is negative before the removal and positive after, although it is estimated from fewer than 30 laps and is not significant in either case. Even after removing outliers, the degradation rate of C3 is slightly lower than that of the harder C2, which motivates the joint model of Appendix~\ref{app:regTC_tirecoef_matrix_completion}.

\begin{table}[htbp]
\centering
\small
\begin{tabular}{l r@{}l r r@{}l r}
\toprule
 & \multicolumn{6}{c}{Model} \\
 \cmidrule(lr){2-7}
 & \multicolumn{3}{c}{With outliers} & \multicolumn{3}{c}{Without outliers} \\
 \cmidrule(lr){2-4}\cmidrule(lr){5-7}
Variable & \multicolumn{2}{c}{Coef.} & Std. Err. &
          \multicolumn{2}{c}{Coef.} & Std. Err. \\
\midrule
Intercept & 102.3083 & *** & (0.191) & 101.8442 & *** & (0.089) \\
DRS & -0.3834 & *** & (0.137) & -0.3765 & *** & (0.064) \\
FirstLap\_pos & -0.1156 & *** & (0.030) & 0.3717 & *** & (0.061) \\
Interval\_behind & 0.4954 & *** & (0.172) & 0.6740 & *** & (0.081) \\
Interval\_front & 1.8054 & *** & (0.274) & 1.9258 & *** & (0.129) \\
LapNumber & -0.0913 & *** & (0.009) & -0.0664 & *** & (0.004) \\
LapNumber$^2$ & 0.0006 & *** & (0.000) & 0.0001 &  & (0.000) \\
PitIn & 1.6267 & *** & (0.155) & 0.3997 & *** & (0.075) \\
PitOut & 20.4048 & *** & (0.162) & 20.2851 & *** & (0.079) \\
C3 (ref. C2) & -0.2258 & * & (0.122) & -0.1544 & *** & (0.057) \\
C4 (ref. C2) & 1.1898 & * & (0.672) & 0.8043 & ** & (0.332) \\
Year 2022 & 0.3556 & *** & (0.082) & 0.4237 & *** & (0.039) \\
Year 2023 & -0.0637 &  & (0.082) & -0.0685 & * & (0.038) \\
TireLife $\cdot$ C2 & 0.0611 & *** & (0.008) & 0.0857 & *** & (0.004) \\
TireLife $\cdot$ C3 & 0.0602 & *** & (0.009) & 0.0787 & *** & (0.004) \\
TireLife $\cdot$ C4 & -0.1045 &  & (0.123) & 0.0398 &  & (0.060) \\
Albon & 1.9023 & *** & (0.198) & 1.7684 & *** & (0.092) \\
Alonso & 1.3695 & *** & (0.177) & 1.3088 & *** & (0.083) \\
Bottas & 1.3130 & *** & (0.186) & 1.2742 & *** & (0.087) \\
Gasly & 1.9389 & *** & (0.188) & 1.4930 & *** & (0.088) \\
Giovinazzi & 1.9019 & *** & (0.249) & 1.9433 & *** & (0.116) \\
Hamilton & 0.0843 &  & (0.172) & 0.1482 & * & (0.080) \\
Hülkenberg & 1.5815 & *** & (0.242) & 1.7902 & *** & (0.113) \\
Latifi & 2.9836 & *** & (0.199) & 2.7591 & *** & (0.093) \\
Leclerc & 0.6149 & *** & (0.172) & 0.5937 & *** & (0.080) \\
Magnussen & 2.1983 & *** & (0.197) & 2.2395 & *** & (0.092) \\
Mazepin & 4.6397 & *** & (0.250) & 4.3056 & *** & (0.120) \\
Norris & 0.9059 & *** & (0.173) & 0.9570 & *** & (0.080) \\
Ocon & 2.7198 & *** & (0.204) & 2.1169 & *** & (0.097) \\
Pérez & 0.4450 & ** & (0.173) & 0.4811 & *** & (0.080) \\
Piastri & 3.5756 & *** & (0.500) & 1.0198 & *** & (0.246) \\
Räikkönen & 1.9776 & *** & (0.249) & 1.8211 & *** & (0.117) \\
Ricciardo & 1.7759 & *** & (0.174) & 1.7124 & *** & (0.081) \\
Russell & 1.2953 & *** & (0.173) & 1.2940 & *** & (0.080) \\
Sainz & 0.6347 & *** & (0.194) & 0.6588 & *** & (0.090) \\
Sargeant & 1.7842 & *** & (0.244) & 1.6352 & *** & (0.114) \\
Schumacher & 2.7114 & *** & (0.199) & 2.6256 & *** & (0.093) \\
Stroll & 1.5572 & *** & (0.184) & 1.6824 & *** & (0.086) \\
Tsunoda & 1.6039 & *** & (0.174) & 1.6050 & *** & (0.081) \\
Vettel & 1.8525 & *** & (0.198) & 1.7463 & *** & (0.092) \\
Zhou & 2.1083 & *** & (0.199) & 1.9340 & *** & (0.093) \\ \midrule
Observations & \multicolumn{3}{c}{2812} & \multicolumn{3}{c}{2764} \\
$R^2$        & \multicolumn{3}{c}{0.877} & \multicolumn{3}{c}{0.966} \\
\bottomrule
\end{tabular}
\caption{Regression coefficients under clear-track conditions for the 2024 United States Grand Prix, estimated with 2021--2023 data. Left: estimated with all data points; Right: estimated after removing outliers. The reference driver is Verstappen. Notes: *** $p<0.01$, ** $p<0.05$, * $p<0.1$.}
\label{tab:regression_cleanlap}
\end{table}

\FloatBarrier
\section{Outlier lap times of the United States GP}\label{app:residuals}

Figure~\ref{fig:laptimes_outliers} shows the clear-track lap times of the 2021--2023 United States GP used to estimate the 2024 coefficients, where outliers are marked with a different color. Most outliers are lap times considerably higher than those of other drivers, and many of them occur in the first lap, when cars are still very close. These laps include driver errors, minor collisions and other incidents that add time. Unusually fast laps do not deviate enough from the rest to be classified as outliers. Laps in which a car enters or leaves the pit lane are also much slower, but they are not outliers, since Equation~\eqref{eq:reg} accounts for them through PitIn and PitOut.

\begin{figure}[htbp]
\centering
\includegraphics[width=\linewidth]{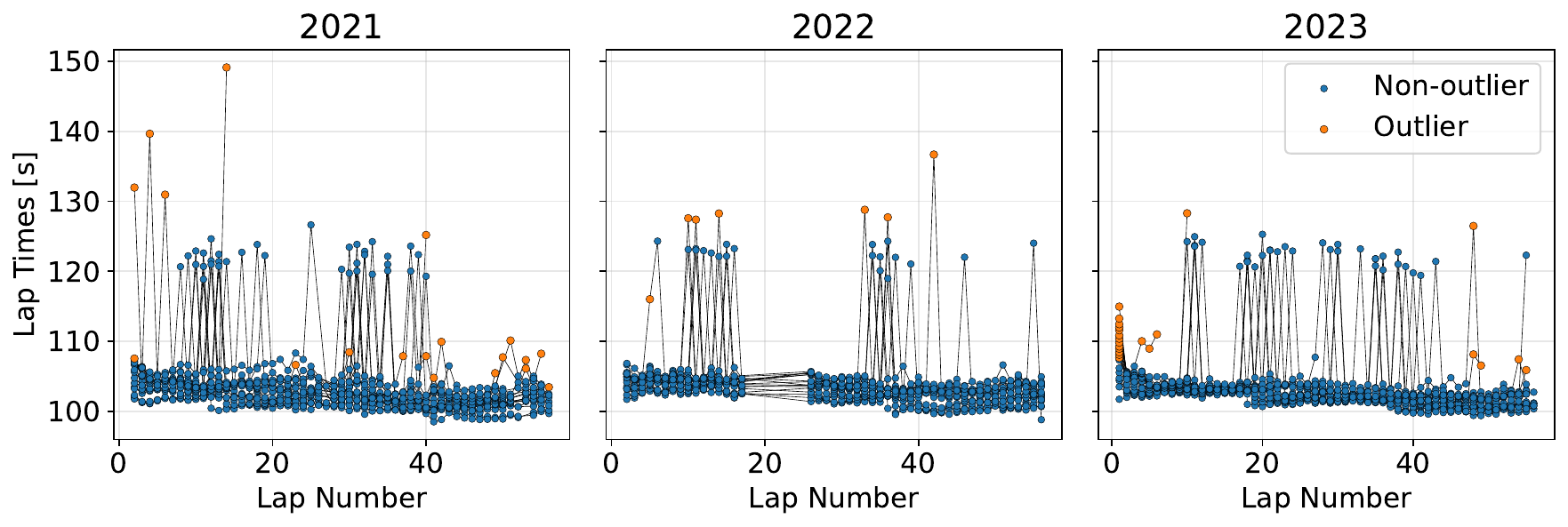}
\caption{Clear-track lap times of the United States GP in 2021, 2022 and 2023. Laps removed by Algorithm~\ref{alg:IQR} are shown in orange.}
\label{fig:laptimes_outliers}
\end{figure}}

\FloatBarrier
\section{Correction of drivers' estimated coefficients for Belgian GP 2024}\label{app:updatebelgium}

\new{Table~\ref{tab:regression_pitstops} reports the coefficients of Equation~\eqref{eq:regcoef} estimated with the races prior to the 2024 Belgian GP. Since its qualifying session was run in wet conditions, the specification without $X_{\text{Q}}$ is used for this race; the full specification is shown for reference. The three previous races and the free practice sessions carry most of the weight, and all coefficients are statistically significant.

\begin{table}[htbp]
\centering
\begin{tabular}{l r@{}l r r@{}l r}
\toprule
 & \multicolumn{3}{c}{Full specification} & \multicolumn{3}{c}{Without $X_{\text{Q}}$ (used)} \\ \cmidrule(lr){2-4}\cmidrule(lr){5-7}
Variable & \multicolumn{2}{c}{Coef.} & Std. Err. & \multicolumn{2}{c}{Coef.} & Std. Err. \\ \midrule
Intercept & 0.0018 & *** & (0.0004) & 0.0016 & *** & (0.0004) \\
$X_{\text{Prev}_1}$ & 0.202 & *** & (0.027) & 0.213 & *** & (0.029) \\
$X_{\text{Prev}_2}$ & 0.213 & *** & (0.027) & 0.201 & *** & (0.028) \\
$X_{\text{Prev}_3}$ & 0.202 & *** & (0.028) & 0.224 & *** & (0.029) \\
$X_{\text{FP}}$ & 0.275 & *** & (0.027) & 0.271 & *** & (0.028) \\
$X_{\text{Q}}$ & 0.039 & *** & (0.009) & \multicolumn{3}{c}{} \\ \midrule
Observations & \multicolumn{3}{c}{929} & \multicolumn{3}{c}{1022} \\
$R^2$ & \multicolumn{3}{c}{0.578} & \multicolumn{3}{c}{0.516} \\ \bottomrule
\end{tabular}
\caption{Coefficients of Equation~\eqref{eq:regcoef} estimated with the races prior to the 2024 Belgian GP. The specification without $X_{\text{Q}}$ is the one used for this race. Notes: *** $p<0.01$, ** $p<0.05$, * $p<0.1$.}
\label{tab:regression_pitstops}
\end{table}}

\FloatBarrier
\section{Overtaking estimated coefficients}\label{app:overtake}

\rev{The overtaking logit is estimated per Grand Prix on 2021--2023 pairs of consecutive cars with a predicted gap between $-30$ and $5$~s, excluding laps 1 to 4 and laps with SC, VSC or red flag. In the simulation the logit decides every battle under clear-track conditions, whatever the gap; the fitted probability $\tau$ is used directly as the success probability of an overtaking attempt (Algorithm~\ref{alg:simevolution}); no decision threshold is applied. The last column of Table~\ref{tab:overtake_logit} reports the empirical overtaking rate in the estimation sample.}

\begin{table}[htbp]
\centering
\begin{tabular}{l r@{}l r@{}l c c c}
\toprule
 & \multicolumn{4}{c}{Estimated Coefficients} & \multicolumn{2}{c}{Training data [observations]} & \multirow{2}{*}{\shortstack{Overtake \\ \rev{rate}}} \\
 \cmidrule(lr){2-5}\cmidrule(lr){6-7}
Grand Prix    & \multicolumn{2}{c}{$\gamma_{0}$} & \multicolumn{2}{c}{$\gamma_{\text{PredGap}}$} & Overtakes & Total & \\ \midrule

Bahrain & -1.1439 & *** & -1.8735 & *** & 235 & 1475 & 0.159 \\
Belgian & -0.7949 & *** & -2.6221 & *** & 183 & 964 & 0.190 \\
Dutch & -1.9121 & *** & -1.8383 & *** & 130 & 1551 & 0.084 \\
Hungarian & -1.6567 & *** & -2.4366 & *** & 206 & 1969 & 0.105 \\
Miami & -0.8381 & *** & -2.4429 & *** & 133 & 1302 & 0.102 \\
Saudi Arabian & -1.3699 & *** & -1.5285 & *** & 110 & 1377 & 0.080 \\
Singapore & -1.9813 & *** & -2.0974 & *** & 44 & 844 & 0.052 \\
United States & -0.6854 & *** & -2.2913 & *** & 199 & 1641 & 0.121 \\

\bottomrule
\end{tabular}
\caption{Overtake logistic regression coefficients for selected GPs, laps with overtake in training data, and \rev{empirical overtaking rate in the estimation sample}. Notes: *** $p<0.01$, ** $p<0.05$, * $p<0.1$.}
\label{tab:overtake_logit}
\end{table}

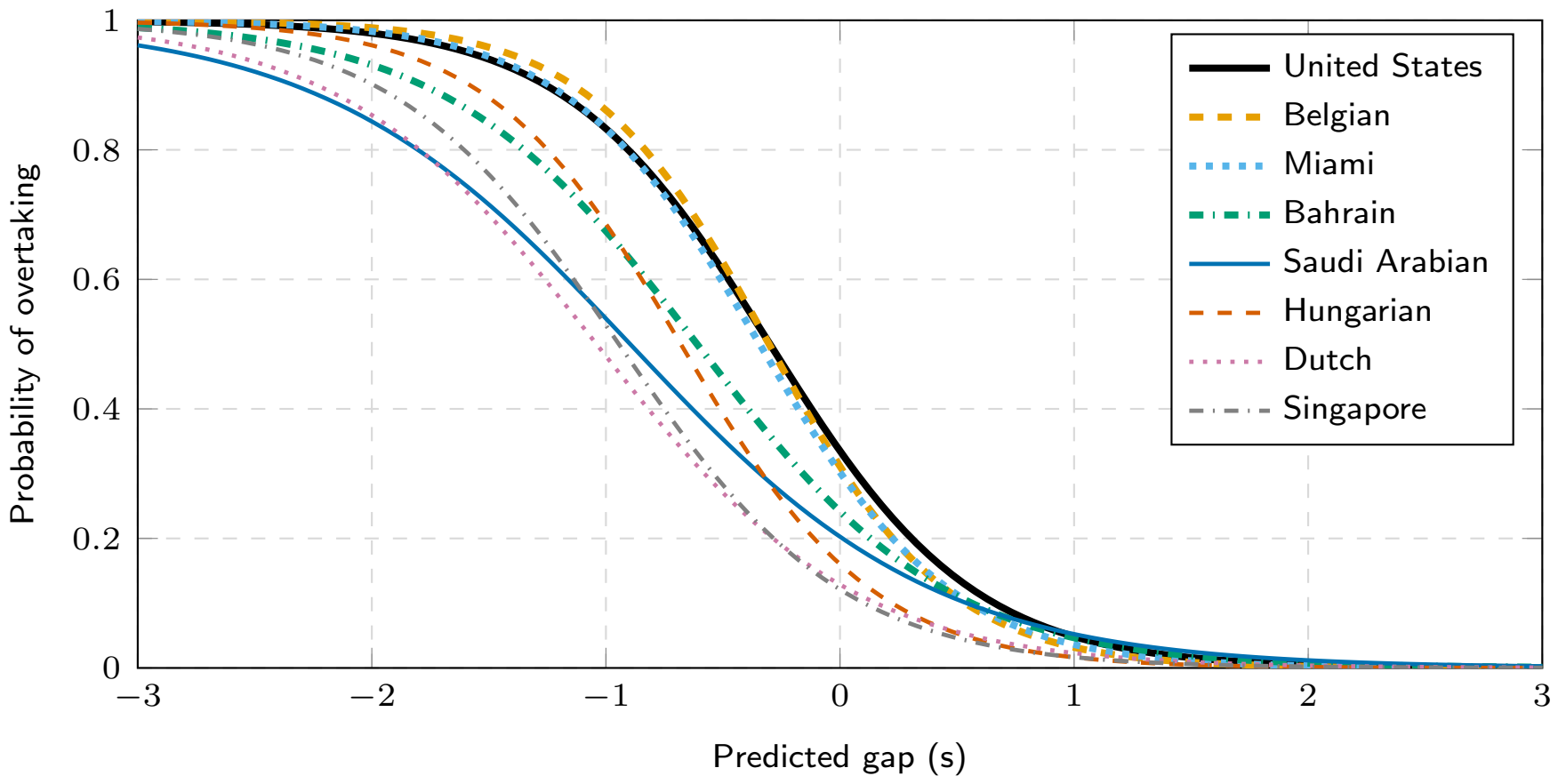
\begin{figure}[htbp]
\centering
\definecolor{oiOrange}{HTML}{E69F00}\definecolor{oiSky}{HTML}{56B4E9}\definecolor{oiGreen}{HTML}{009E73}
\definecolor{oiBlue}{HTML}{0072B2}\definecolor{oiVerm}{HTML}{D55E00}\definecolor{oiPurple}{HTML}{CC79A7}
\begin{tikzpicture}
\begin{axis}[
    width=0.92\linewidth, height=6.4cm,
    xlabel={Predicted gap (s)}, ylabel={Probability of overtaking},
    xmin=-3, xmax=3, ymin=0, ymax=1,
    domain=-3:3, samples=120, no markers,
    ymajorgrids=true, xmajorgrids=true, grid style={dashed, gray!30},
    tick label style={font=\footnotesize\sffamily}, label style={font=\footnotesize\sffamily},
    legend style={at={(0.98,0.98)}, anchor=north east, font=\footnotesize\sffamily, cells={anchor=west}, legend image post style={xscale=1.0}},
]
\addplot[black, line width=1.5pt] {1/(1+exp(-(-0.6854-2.2913*x)))};
\addlegendentry{United States}
\addplot[oiOrange, dashed, line width=1.5pt] {1/(1+exp(-(-0.7949-2.6221*x)))};
\addlegendentry{Belgian}
\addplot[oiSky, dotted, line width=1.5pt] {1/(1+exp(-(-0.8381-2.4429*x)))};
\addlegendentry{Miami}
\addplot[oiGreen, dashdotted, line width=1.5pt] {1/(1+exp(-(-1.1439-1.8735*x)))};
\addlegendentry{Bahrain}
\addplot[oiBlue, line width=0.9pt] {1/(1+exp(-(-1.3699-1.5285*x)))};
\addlegendentry{Saudi Arabian}
\addplot[oiVerm, dashed, line width=0.9pt] {1/(1+exp(-(-1.6567-2.4366*x)))};
\addlegendentry{Hungarian}
\addplot[oiPurple, dotted, line width=0.9pt] {1/(1+exp(-(-1.9121-1.8383*x)))};
\addlegendentry{Dutch}
\addplot[gray, dashdotted, line width=0.9pt] {1/(1+exp(-(-1.9813-2.0974*x)))};
\addlegendentry{Singapore}
\end{axis}
\end{tikzpicture}
\caption{Probability of overtaking as a function of the predicted gap, from the coefficients of Table~\ref{tab:overtake_logit}. Circuits are listed in decreasing order of the probability at a gap of zero.}
\label{fig:overtake_prob}
\end{figure}

Figure~\ref{fig:overtake_prob} shows the probability of overtaking implied by the coefficients of Table~\ref{tab:overtake_logit} as a function of the predicted gap $X_{\text{PredGap}}$, where a negative gap means that the trailing car is predicted to finish the lap ahead. At a predicted gap of zero, the probability ranges from $12\%$ (Singapore) to $34\%$ (United States). The trailing car must be predicted to finish ahead by between $0.3$ and $1.0$~s for the probability to reach $50\%$, and by about $3$~s for it to exceed $95\%$ at every circuit.

\FloatBarrier
\section{Rivals' Pit-Stop Decision-Making}\label{app:rivals_strategy}
Table~\ref{tab:logit_pitin} reports the estimated coefficients of the logistic regression in Equation~\eqref{eq:logit_rivals_stop_continue}. \new{Laps under a safety car (SC) or virtual safety car (VSC) and a worse track position are all associated with a higher likelihood of making a pit stop, whereas stops are less likely in the first ten laps and much less likely in the last ten.} If a driver has already used two different tire compounds, the odds of making a stop decrease. In general, the softer the tire compound, the more likely a driver is to stop earlier. \rev{In the simulation, the logit is evaluated only on laps in which a pit stop is admissible, i.e., clear-track laps and the first lap of a SC \new{or VSC} period (Algorithm~\ref{alg:simevolution}). The model is estimated with one observation per driver and SC or VSC period, and stops are identified by a change of tire set, since the timing data also flag laps in which the whole field was led through the pit lane behind the SC.}

\begin{table}[htbp]
\centering
\begin{tabular}{l r@{}l r}
\toprule
Variable & \multicolumn{2}{c}{Coef.} & \multicolumn{1}{c}{Std. Err.} \\
\midrule
Intercept & -5.6162 & *** & (0.152) \\
Bahrain GP & 1.5180 & *** & (0.149) \\
Belgian GP & 1.6129 & *** & (0.170) \\
Dutch GP & 0.4334 & *** & (0.150) \\
Hungarian GP & 0.0325 &  & (0.154) \\
Miami GP & -0.7288 & *** & (0.214) \\
Saudi Arabian GP & 0.1041 &  & (0.198) \\
Singapore GP & -1.4500 & *** & (0.257) \\
United States GP & 1.3122 & *** & (0.155) \\
SC & 3.5309 & *** & (0.115) \\
VSC & 3.1374 & *** & (0.144) \\
Used two compounds & -1.0413 & *** & (0.065) \\
Early laps & -0.4136 & *** & (0.097) \\
Late laps & -2.0614 & *** & (0.115) \\
TireLife $\cdot$ C1 & 0.1380 & *** & (0.009) \\
TireLife $\cdot$ C2 & 0.1180 & *** & (0.005) \\
TireLife $\cdot$ C3 & 0.1461 & *** & (0.004) \\
TireLife $\cdot$ C4 & 0.1893 & *** & (0.006) \\
TireLife $\cdot$ C5 & 0.3002 & *** & (0.014) \\
Position & 0.0490 & *** & (0.005) \\
\midrule
Observations & \multicolumn{3}{c}{60,433} \\
Pseudo $R^2$ & \multicolumn{3}{c}{0.216} \\
\bottomrule
\end{tabular}
\caption{Logit regression results for pit stop decision. Each SC or VSC period contributes one observation per driver, on its first lap. Notes: *** $p<0.01$, ** $p<0.05$, * $p<0.1$. Reference categories and the indicators of the circuits not considered omitted.}
\label{tab:logit_pitin}
\end{table}

Table~\ref{tab:logit_compound_choice} shows the estimated coefficients of the conditional logistic regression in Equation~\eqref{eq:cond_logit_rivals}. The negative sign of the last coefficient indicates that when the driver has used one compound only, it is more likely that a different compound is selected, all else equal.
\begin{table}[htbp]
\centering
\begin{tabular}{l r@{}l r}
\toprule
Variable & \multicolumn{2}{c}{Coef.} & \multicolumn{1}{c}{Std. Err.} \\
\midrule

Alt: C0 & -1.3493 & ** & (0.540) \\
Alt: C1 & -0.0278 &  & (0.610) \\
Alt: C2 & -1.0243 & ** & (0.406) \\
Alt: C4 & 1.2418 & *** & (0.443) \\
Alt: C5 & 3.4327 & *** & (1.044) \\

Alt: C0 $\cdot$ LapsLeft & 0.0265 & * & (0.015) \\
Alt: C1 $\cdot$ LapsLeft & 0.0054 &  & (0.019) \\
Alt: C2 $\cdot$ LapsLeft & 0.0278 & * & (0.015) \\
Alt: C4 $\cdot$ LapsLeft & -0.0793 & *** & (0.017) \\
Alt: C5 $\cdot$ LapsLeft & -0.2823 & *** & (0.063) \\

Alt: Same \& only one comp. used & -1.5111 & *** & (0.296) \\

\midrule
Observations & \multicolumn{3}{c}{1,742} \\
Pseudo $R^2$ & \multicolumn{3}{c}{0.311} \\
\bottomrule
\end{tabular}
\caption{Conditional logit results for tire compound choice conditional on pit stop. Notes: *** $p<0.01$, ** $p<0.05$, * $p<0.1$. \new{Base alternative (C3) omitted. Clustered standard errors at race level (63 races).}}
\label{tab:logit_compound_choice}
\end{table}

\FloatBarrier
\section{DRL Performance in Simplified Deterministic Instances}\label{app:opt_gap_simple_instance}

Table~\ref{tab:opt_gap} reports, for the simplified setting with one driver and no uncertainty, how much slower the full-race strategy of each DRL agent is than the optimal strategy of the QMIP. The gap is in seconds, and a zero means that the agent reaches the optimal race time. Since different strategies can have the same race time, the table complements the stints shown in Figure~\ref{fig:stints_alone}. The first eight rows give, for each circuit, the seed with the smallest gap among the five trained that satisfy the two-compound rule, and the row ``Mean'' averages them over circuits. The last row gives the median over the five seeds, averaged over circuits. It is larger than the mean of the selected seeds for every method, so the result depends on the training run.

\begin{table}[htbp]
\centering
\begin{tabular}{lccccc}
\toprule
Grand Prix & Double DQN & A2C & TRPO & PPO & Recurrent PPO \\
\midrule
Bahrain & 0.000 & 0.185 & 4.919 & 0.068 & 0.000 \\
Belgian & 0.057 & 0.057 & 1.217 & 0.000 & 0.000 \\
Dutch & 0.000 & 1.231 & 1.408 & 1.231 & 0.000 \\
Hungarian & 0.000 & 0.128 & 5.374 & 1.217 & 0.000 \\
Miami & 0.421 & 0.889 & 3.438 & 0.885 & 0.128 \\
Saudi Arabian & 0.000 & 3.471 & 5.068 & 0.000 & 0.000 \\
Singapore & 0.000 & 4.971 & 5.034 & 0.012 & 0.000 \\
United States & 0.102 & 0.082 & 2.237 & 0.279 & 0.000 \\
\midrule
Mean & 0.072 & 1.377 & 3.587 & 0.461 & 0.016 \\
\midrule
Median over seeds & 0.596 & 2.214 & 6.893 & 2.503 & 0.903 \\
\bottomrule
\end{tabular}
\caption{\new{Time difference (in seconds) between each DRL model strategy and the global optimal race strategy computed via the QMIP model. Positive values indicate how much slower the agent's full-race strategy is compared to the QMIP optimum. Each entry corresponds to the seed with the smallest gap among five. \textit{Median over seeds}: median gap across the five seeds of each circuit, averaged over circuits.}}\label{tab:opt_gap}
\end{table}

\FloatBarrier
\section{Drivers Who Did Not Finish (DNF) per Race}\label{app:DNF}

Table~\ref{tab:DNFs_2024} lists the drivers who did not finish each of the eight 2024 races. In the simulations, rivals that retired in the real race retire on the same lap (Section~\ref{sec:agent_competition}). The second column counts all the retirements of the race, and the third counts those of drivers who started ahead of the agent. In these races, every retirement is of a driver who started behind the agent.
\begin{table}[htbp]
\centering
\begin{tabular}{lcc}
\toprule
        &  \multicolumn{2}{c}{DNF}  \\ \cmidrule{2-3}
Grand Prix     &  Total   & Ahead of Agent   \\ \midrule
Bahrain        &     0        &   0         \\
Belgian        &     1        &   0         \\
Dutch          &     0        &   0         \\
Hungarian      &     1        &   0         \\
Miami          &     1        &   0         \\
Saudi Arabian  &     2        &   0         \\
Singapore      &     2        &   0         \\
United States  &     1        &   0         \\
\bottomrule
\end{tabular}
\caption{Number of drivers who did not finish (DNF) per race in the 2024 season. DNF Ahead of Agent: number of DNFs among drivers who started ahead of the agent's grid position.}
\label{tab:DNFs_2024}
\end{table}

\FloatBarrier
\new{\section{Mean Final Positions with Confidence Intervals}\label{app:ci_2024}

Table~\ref{tab:ci_2024} reports the mean final position of each method in the eight circuits, in the fully simulated races (top) and in the races with the real realizations (bottom). A lower position is better. The row ``Average'' is the mean over the circuits, and ``Difference with QMIP'' is the mean paired difference averaged over the circuits, so a negative value means that the method finishes ahead of the QMIP. The marks are not corrected for multiple comparisons. The intervals do not include the variability across training seeds, so the differences with the QMIP refer to the selected agent of each method.
\begin{table}[htbp]
\centering
\resizebox{\textwidth}{!}{%
\begin{tabular}{lcccccc}
\toprule
GP & Double DQN & A2C & TRPO & PPO & Rec. PPO & QMIP \\
\midrule
\multicolumn{7}{l}{\textit{Fully simulated races (Figure~\ref{fig:stats_drivers_stochastic_2024})}} \\
\midrule
Bahrain & 5.34 $\pm$ 0.24$^{\dagger}$ & 4.69 $\pm$ 0.23$^{\dagger}$ & 4.91 $\pm$ 0.24$^{\dagger}$ & 4.04 $\pm$ 0.19$^{*}$ & 4.13 $\pm$ 0.19$^{*}$ & 4.34 $\pm$ 0.20 \\
Belgian & 10.74 $\pm$ 0.25$^{*}$ & 10.71 $\pm$ 0.24$^{*}$ & 11.79 $\pm$ 0.29$^{\dagger}$ & 10.77 $\pm$ 0.25$^{*}$ & 10.80 $\pm$ 0.25$^{*}$ & 11.18 $\pm$ 0.25 \\
Dutch & 7.67 $\pm$ 0.26 & 7.98 $\pm$ 0.26$^{\dagger}$ & 7.83 $\pm$ 0.30 & 6.79 $\pm$ 0.22$^{*}$ & 6.70 $\pm$ 0.22$^{*}$ & 7.72 $\pm$ 0.25 \\
Hungarian & 8.47 $\pm$ 0.25$^{\dagger}$ & 8.37 $\pm$ 0.24$^{\dagger}$ & 8.28 $\pm$ 0.24$^{\dagger}$ & 7.73 $\pm$ 0.22$^{*}$ & 7.84 $\pm$ 0.23 & 8.02 $\pm$ 0.25 \\
Miami & 12.98 $\pm$ 0.29 & 13.85 $\pm$ 0.29$^{\dagger}$ & 14.50 $\pm$ 0.28$^{\dagger}$ & 12.76 $\pm$ 0.28$^{*}$ & 13.06 $\pm$ 0.28 & 12.99 $\pm$ 0.27 \\
Saudi Arabian & 4.61 $\pm$ 0.18 & 4.45 $\pm$ 0.19 & 4.71 $\pm$ 0.17$^{\dagger}$ & 4.41 $\pm$ 0.16 & 4.26 $\pm$ 0.16$^{*}$ & 4.51 $\pm$ 0.16 \\
Singapore & 6.65 $\pm$ 0.22$^{*}$ & 6.74 $\pm$ 0.22 & 8.62 $\pm$ 0.33$^{\dagger}$ & 6.42 $\pm$ 0.21$^{*}$ & 6.40 $\pm$ 0.21$^{*}$ & 6.93 $\pm$ 0.21 \\
United States & 9.20 $\pm$ 0.23$^{\dagger}$ & 9.07 $\pm$ 0.23$^{\dagger}$ & 9.06 $\pm$ 0.25$^{\dagger}$ & 8.89 $\pm$ 0.23 & 9.00 $\pm$ 0.23$^{\dagger}$ & 8.74 $\pm$ 0.23 \\
\midrule
Average & 8.21 $\pm$ 0.09$^{\dagger}$ & 8.23 $\pm$ 0.08$^{\dagger}$ & 8.71 $\pm$ 0.09$^{\dagger}$ & 7.73 $\pm$ 0.08$^{*}$ & 7.77 $\pm$ 0.08$^{*}$ & 8.05 $\pm$ 0.08 \\
Difference with QMIP & $+$0.16 $\pm$ 0.08 & $+$0.18 $\pm$ 0.08 & $+$0.66 $\pm$ 0.09 & $-$0.33 $\pm$ 0.07 & $-$0.28 $\pm$ 0.07 & -- \\
\midrule
\multicolumn{7}{l}{\textit{Real race realizations (Figure~\ref{fig:boxplot_drivers_2024})}} \\
\midrule
Bahrain & 5.77 $\pm$ 0.22$^{\dagger}$ & 4.60 $\pm$ 0.22$^{*}$ & 4.95 $\pm$ 0.23 & 4.77 $\pm$ 0.22$^{*}$ & 4.83 $\pm$ 0.23$^{*}$ & 5.00 $\pm$ 0.23 \\
Belgian & 12.19 $\pm$ 0.27 & 12.19 $\pm$ 0.23 & 12.37 $\pm$ 0.27 & 12.30 $\pm$ 0.25 & 12.89 $\pm$ 0.27$^{\dagger}$ & 12.31 $\pm$ 0.25 \\
Dutch & 11.06 $\pm$ 0.24$^{\dagger}$ & 9.59 $\pm$ 0.21$^{*}$ & 9.65 $\pm$ 0.21$^{*}$ & 9.61 $\pm$ 0.20$^{*}$ & 9.83 $\pm$ 0.21$^{*}$ & 10.73 $\pm$ 0.22 \\
Hungarian & 8.55 $\pm$ 0.18$^{*}$ & 8.46 $\pm$ 0.16$^{*}$ & 8.45 $\pm$ 0.16$^{*}$ & 8.45 $\pm$ 0.15$^{*}$ & 9.54 $\pm$ 0.25$^{\dagger}$ & 9.35 $\pm$ 0.23 \\
Miami & 10.58 $\pm$ 0.30$^{*}$ & 11.01 $\pm$ 0.35$^{*}$ & 12.62 $\pm$ 0.35$^{\dagger}$ & 11.07 $\pm$ 0.33$^{*}$ & 12.57 $\pm$ 0.29$^{\dagger}$ & 12.27 $\pm$ 0.28 \\
Saudi Arabian & 7.20 $\pm$ 0.20$^{\dagger}$ & 5.34 $\pm$ 0.14 & 5.66 $\pm$ 0.18$^{\dagger}$ & 5.25 $\pm$ 0.14 & 5.25 $\pm$ 0.14 & 5.26 $\pm$ 0.14 \\
Singapore & 9.30 $\pm$ 0.16$^{\dagger}$ & 8.31 $\pm$ 0.18$^{\dagger}$ & 9.04 $\pm$ 0.20$^{\dagger}$ & 8.41 $\pm$ 0.17$^{\dagger}$ & 7.92 $\pm$ 0.16 & 7.96 $\pm$ 0.17 \\
United States & 9.74 $\pm$ 0.21$^{*}$ & 10.13 $\pm$ 0.21 & 10.78 $\pm$ 0.25$^{\dagger}$ & 9.39 $\pm$ 0.19$^{*}$ & 9.49 $\pm$ 0.20$^{*}$ & 10.27 $\pm$ 0.22 \\
\midrule
Average & 9.30 $\pm$ 0.08$^{\dagger}$ & 8.70 $\pm$ 0.08$^{*}$ & 9.19 $\pm$ 0.08 & 8.66 $\pm$ 0.08$^{*}$ & 9.04 $\pm$ 0.08$^{*}$ & 9.14 $\pm$ 0.08 \\
Difference with QMIP & $+$0.16 $\pm$ 0.06 & $-$0.44 $\pm$ 0.06 & $+$0.05 $\pm$ 0.07 & $-$0.49 $\pm$ 0.06 & $-$0.10 $\pm$ 0.05 & -- \\
\bottomrule
\end{tabular}}
\caption{Mean final position by solution method and circuit, 2024 season, over $500$ simulations per circuit, with $95\%$ confidence intervals. Since all methods are evaluated on the same $500$ random seeds, the difference with the QMIP is computed race by race. The intervals reflect the variability across simulated races for the trained agent of each method, and not the variability across training seeds or circuits. $^{*}$/$^{\dagger}$~Better/worse than the QMIP at the 5\% level (paired $t$-test, i.e., the 95\% confidence interval of the paired difference excludes zero).}\label{tab:ci_2024}
\end{table}}

\FloatBarrier
\section{Tire Stints of the Candidate Actions}\label{app:decision_stints}

\begin{figure}[htbp]
\centering
\includegraphics[width=1.0\linewidth]{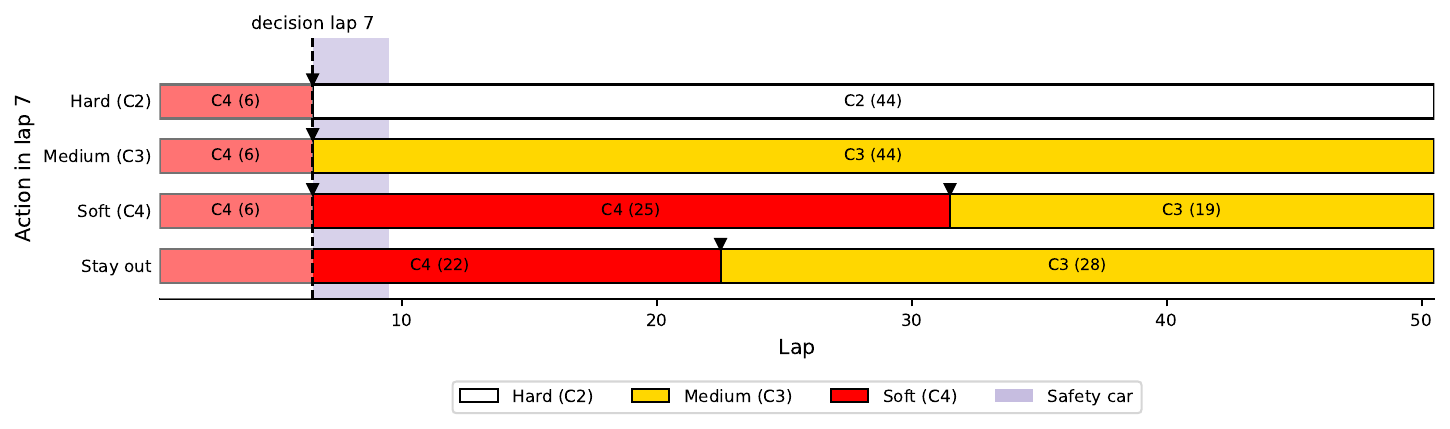}
\caption{Tire stints implied by each candidate action at lap $7$ of the 2024 Saudi Arabian Grand Prix, with the PPO policy driving in place of Lando Norris. Each bar shows one representative continuation under the PPO policy, with the same realization of race events across actions. Numbers in parentheses give stint lengths in laps.}
\label{fig:decision_stints}
\end{figure}

Figure~\ref{fig:decision_stints} shows the tire stints that follow each of the four candidate actions at the decision point of Figure~\ref{fig:decision_support}. In every bar, the first segment is the stint on soft tires run before the decision, the dashed line marks the decision lap, the shaded band is the safety car period, and the triangles mark the pit stops. Fitting hard or medium tires at lap $7$ is enough to finish the race, since the new set runs the remaining $44$ laps. Fitting soft tires again requires a second stop, to medium tires, after a stint of $25$ laps. If the driver stays out, the soft stint ends after $22$ laps in a stop to medium tires. These second stops are why the two actions lose points in Figure~\ref{fig:decision_support}. Each bar is one continuation under the PPO policy, with the same race events in the four of them, so the bars show strategies, not average outcomes.

\begin{figure}[htbp]
\centering
\includegraphics[width=1.0\linewidth]{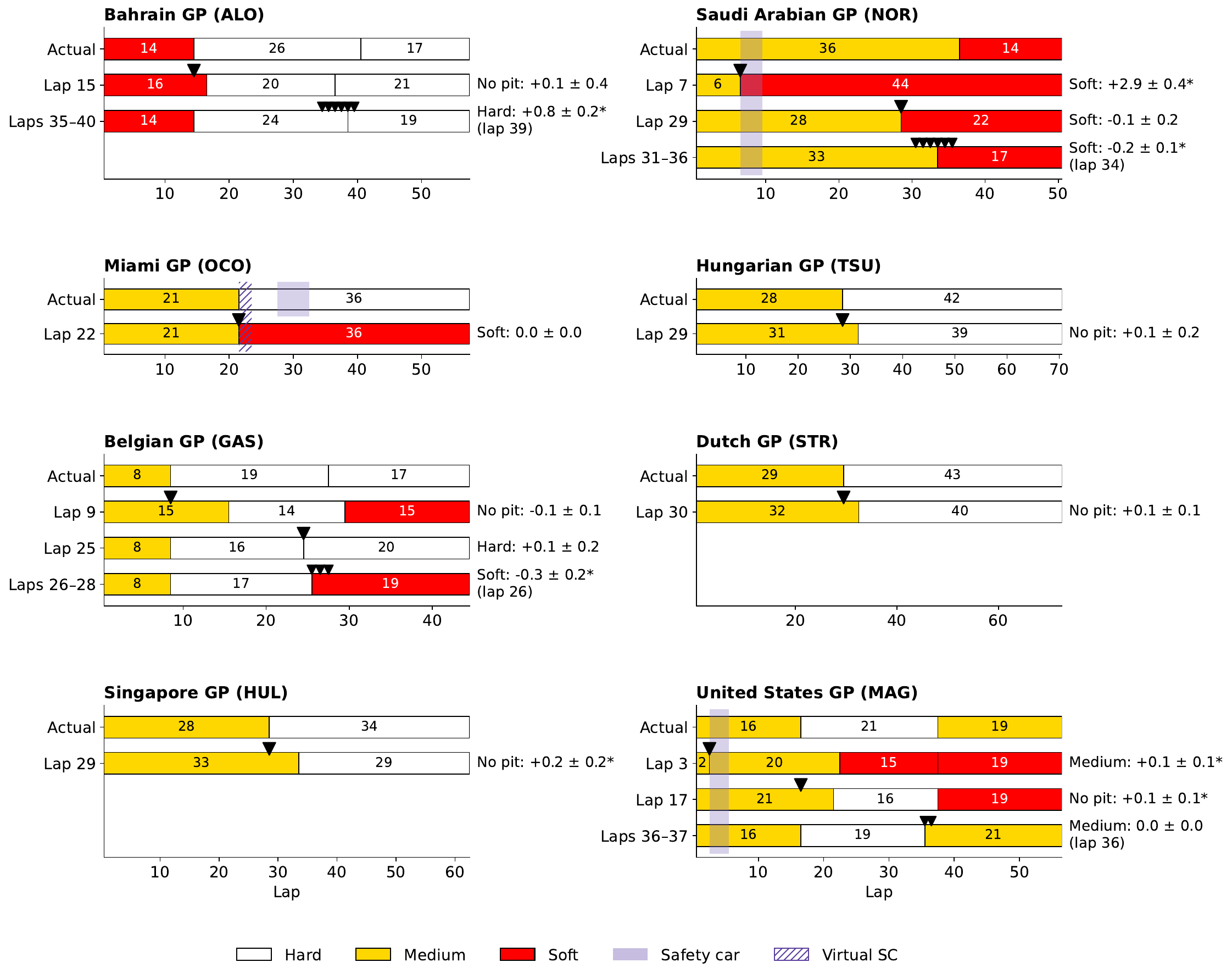}
\caption{Actual 2024 tire stints of each driver and stints implied by the DRL recommendation at every lap where it differs from the actual decision. Each lower bar replays the race up to the lap marked by the triangle and shows one representative continuation under the PPO policy, without further yellow flags. Consecutive laps with the same recommendation share a bar, drawn from the lap with the largest difference. The annotation on the right reports the recommended action and the mean paired difference in points against the actual decision at that lap (given in parentheses for grouped bars), over $M = 500$ rollouts with simulated yellow flags, $\pm 2$ SE; positive values favor the DRL recommendation and negative values favor the driver's actual decision, and $^{*}$ marks differences larger than 2 SE. The differences come from the rollouts used to select the best action in Figure~\ref{fig:realhist_best}, so they differ slightly from those reported there.}
\label{fig:realhist}
\end{figure}

\FloatBarrier
\new{\section{Review of the 2024 Decisions: Further Training and Results for All Races}\label{app:realhist}

The policies evaluated in Section~\ref{sec:agent_competition} choose their own starting compound and, in most races, always start on the same one. When the race history fixes a different starting compound, as for Norris in the 2024 Saudi Arabian Grand Prix, the look-ahead rollouts reach states rarely visited during training, in which the policy may propose a second stint on the starting compound before the two-compound rule is met, so that the restriction to feasible actions in Algorithm~\ref{alg:counterfactual_rollouts} has to override it.

For the analysis of Figure~\ref{fig:realhist}, training resumes from the selected PPO policy of each race, including its observation and reward normalization statistics, in the same environment and with the hyperparameters of Table~\ref{tab:params_all_rivals}, except for the following. In half of the episodes, the starting compound is drawn uniformly among the available ones, and it is chosen by the policy otherwise; the same mixture is used in the periodic evaluation episodes. The learning rate is kept constant at $10^{-4}$, and each run lasts at most $1.5$ million steps, with an evaluation every $5{,}000$ steps per environment and the same early stopping as in the original training, keeping the best checkpoint. We train ten seeds per race and select one with the validation rule of Appendix~\ref{app:architecture}, among the runs that also meet the two-compound rule in all of $200$ validation races when forced to start on the driver's actual 2024 compound.

Over the eight races, the further training changes the mean final position in the two scenarios of Figures~\ref{fig:stats_drivers_stochastic_2024} and~\ref{fig:boxplot_drivers_2024} by less than $0.05$ positions; the policy remains ahead of the QMIP and meets the two-compound rule in all the simulated races, compared with up to $0.2\%$ of failures for the original policy.}

Figure~\ref{fig:realhist} shows the stints implied by the recommendation at every lap where it differs from the actual decision in the eight races. In Singapore, Hülkenberg started on medium tires and stopped on lap $29$ for hard tires ($28$ and $34$ laps). At that lap the recommendation is to stay out, with a gain of $0.2 \pm 0.2$ points, just above two standard errors. Under the policy the driver also makes one stop, from medium to hard tires, but five laps later ($33$ and $29$ laps), so the difference is small. In the Saudi Arabian Grand Prix, Norris started on medium tires, stayed out under the safety car on lap $7$, and stopped on lap $37$ for soft tires ($36$ and $14$ laps). At lap $7$ the recommendation is to fit soft tires, which gains $2.9 \pm 0.4$ points over staying out, and the policy runs them to the end ($6$ and $44$ laps). In Figure~\ref{fig:decision_support} the policy recommends medium tires instead, since in that race it started on softs. In the later laps the recommendation is again soft tires, but the differences are $-0.1 \pm 0.2$ (lap $29$) and $-0.2 \pm 0.1$ (laps $31$ to $36$), so they do not support stopping earlier than Norris did.

\FloatBarrier
\section{Computation Times}\label{app:runtimes}

Table~\ref{tab:runtimes} reports the median wall-clock time, in hours, needed to train one policy for one circuit, with at most 5M timesteps. Policies are trained on an AMD EPYC 9224, with four CPU cores in the setting with competition and two without it. Once trained, a PPO decision takes about $0.1$ ms and a QMIP solve at most about $20$ ms, or about $0.2$ s for a whole race, on one core of an Apple M4. A look-ahead recommendation is more demanding, since its cost grows with $M$ and with the remaining laps: at lap $7$ of the case in Figure~\ref{fig:decision_support}, it takes about $50$ s with $M=500$ and $11$ s with $M=100$ in a single process. Since rollouts are independent, they can be distributed across processes with identical results: the look-ahead takes about $13$ s for $M=500$ with eight processes and about $6$ s for $M=100$ with four.

\begin{table}[htbp]
\centering
\begin{tabular}{lccccc}
\toprule
Setting & Double DQN & A2C & TRPO & PPO & Recurrent PPO \\
\midrule
With competition & 1.65 & 1.48 & 1.67 & 1.50 & 1.86 \\
Without competition & 0.37 & 0.16 & 0.28 & 0.25 & 0.83 \\
\bottomrule
\end{tabular}
\caption{Median wall-clock training time per policy and circuit, in hours, by DRL method.}
\label{tab:runtimes}
\end{table}

\FloatBarrier
\end{appendices}

\end{document}